\documentclass[a4paper]{article}

\usepackage{fullpage}
\usepackage[english]{babel}
\usepackage{latexsym,amsfonts,amsmath,amssymb,amsthm,amscd}

\usepackage{hyperref}

\usepackage[dvips]{epsfig}
\usepackage{graphics}

\newcommand{\A}{\mathcal{A}}
\newcommand{\B}{\mathcal{B}}

\newcommand{\Ba}{\mathfrak{B}}
\newcommand{\re}{\mathcal{R}}
\newcommand{\M}{\mathbb{M}}
\newcommand{\U}{\mathbb{U}}
\newcommand{\V}{\mathbb{V}}

\newcommand{\E}{\mathbb{E}_{\mu}}

\newcommand{\p}{\mathcal{P}}

\newcommand{\I}{\mathcal{I}}

\newcommand{\N}{\mathbb{N}}

\newcommand{\Z}{\mathbb{Z}}
\newcommand{\C}{\mathbb{C}}
\newcommand{\az}{\mathcal{A}^{\mathbb{Z}}}
\newcommand{\bz}{\mathcal{B}^{\mathbb{Z}}}

\newcommand{\am}{\mathcal{A}^{\mathbb{M}}}
\newcommand{\an}{\mathcal{A}^{\mathbb{N}}}

\newcommand{\bn}{\mathcal{B}^{\mathbb{N}}}

\newcommand{\dd}{\delta}

\newcommand{\s}{\sigma}
\newcommand{\gs}{\Sigma}
\newcommand{\ds}{D^{\Sigma}}

\newcommand{\F}{F_T}
\newcommand{\Ft}{\widetilde{F_{T}}}

\newcommand{\supp}{\mathrm{supp}}
\newcommand{\Succ}{\mathrm{Succ}}
\newcommand{\Hom}{\mathrm{Hom}}

\newcommand{\im}{\mathrm{Im}}
\newcommand{\Ker}{\mathrm{Ker}}
\newcommand{\id}{\mathrm{Id}}
\newcommand{\Me}{\mathcal{M}}

\newcommand{\na}{\N^{\ast}}

\newcommand{\eg}[1]{\raisebox{-1ex}{${{\displaystyle =}\atop {\scriptstyle{\mathrm{#1}}}}$}}
\newcommand{\rightarrrow}[1]{\raisebox{-1ex}{${{\displaystyle \longrightarrow}\atop {\scriptstyle{\mathrm{#1}}}}$}}
\newcommand{\geeeq}[1]{\raisebox{-1ex}{${{\displaystyle >}\atop {\scriptstyle{\mathrm{#1}}}}$}}
\newcommand{\leeeq}[1]{\raisebox{-1ex}{${{\displaystyle \leq}\atop {\scriptstyle{\mathrm{#1}}}}$}}

\title{Measure rigidity for\\ algebraic bipermutative cellular automata}

\theoremstyle{plain}
\newtheorem{theorem}{Theorem}[section]
\newtheorem{prop}[theorem]{Proposition}
\newtheorem{cor}[theorem]{Corollary}
\newtheorem{lemma}[theorem]{Lemma}

\theoremstyle{definition}

\theoremstyle{remark}
\newtheorem{remark}{Remark}[section]
\newtheorem{ex}[remark]{Example}

\newcounter{claimcount}
\newcommand{\THMfont}[1]{{\sl #1}}
\newcommand{\Claim}[1]{\refstepcounter{claimcount} \vspace{0.3em}               \noindent {\sc Claim \theclaimcount: \ }\THMfont{ #1}}
\newcommand{\bprf}[1][Proof:]{\begin{list}{}    {\setlength{\leftmargin}{0.5em} \setlength{\rightmargin}{0em}  \setlength{\listparindent}{1em}}   \item {\em \hspace{-0.8em}  #1  }}
\newcommand{\eprf}{\end{list}}

\newcommand{\bclaimprf}{\bprf}

\newcommand{\eclaimprf}{ \hfill $\Diamond$~{\scriptsize {\tt Claim~\theclaimcount}}\eprf} 

\date{2005}
\author{Mathieu Sablik \\
  Institut de Math\'ematiques de Luminy \\
  UMR 6206-Campus de Luminy, Case 907\\
  13288 Marseille Cedex 09 \\
   \tt{sablik@iml.univ-mrs.fr}  }

\begin{document}
\maketitle

\begin{abstract}
Let $(\az,F)$ be a bipermutative algebraic cellular automaton. We present conditions which force a probability
measure which is invariant for the $\N\times\Z$-action of $F$ and the shift map $\s$ to be the Haar measure on
$\gs$, a closed shift-invariant subgroup of the Abelian compact group $\az$. This generalizes simultaneously
results of B. Host, A. Maass and S. Mart\'{\i}nez \cite{Host-Maass-Martinez-2003} and M. Pivato
\cite{Pivato-2003}. This result is applied to give conditions which also force an $(F,\s)$-invariant probability
measure to be the uniform Bernoulli measure when $F$ is a particular invertible affine expansive cellular
automaton on~$\an$.
\end{abstract}

\section{Introduction}

Let $F:\am\to\am$ with $\M=\N$ or $\Z$ be a one-dimensional cellular automaton (CA). The study of invariant
measures under the action of $F$ has been addressed from different points of view in the last two decades. As
ergodic theory is the study of invariant measures, it is thus natural to characterize them. In addition, since
$F$ commutes with the shift map $\s$, it is important to describe invariant probability measures for the
semi-group action generated by $F$ and $\s$. We remark that it easy to prove the existence of such measures by
considering a cluster point of the Ces\`aro mean under iteration of $F$ of a $\s$-invariant measure. This
problem is related to Furstenberg's conjecture \cite{Furstenberg-1967} that the Lebesgue measure on the torus is
the unique invariant measure under multiplication by two relatively prime integers. In the algebraic setting,
the study of invariant measures under a group action on a zero-dimensional group like Ledrappier's example
\cite{Ledrappier-1978}, has been extensively considered in \cite{Schmidt-1995} and \cite{Einsiedler-2003}.

The uniform Bernoulli measure has an important role in the study of $(F,\s)$-invariant measures. G.A. Hedlund
has shown in \cite{Hedlund-1969} that a CA is surjective iff the uniform Bernoulli measure on $\am$ is
$(F,\s)$-invariant. Later, D. Lind \cite{Lind-1984} shows for the radius $1$ $\textrm{mod}\, 2$ automaton that
starting from any Bernoulli measure the Ces\`aro mean of the iterates by the CA converges to the uniform
measure. This result is generalized for a large class of algebraic CA and a large class of measures with tools
from stochastic processes in \cite{Maass-Martinez-1998} and \cite{Ferrari-Maass-Martinez-Ney-2000}, and with
harmonic analysis tools in \cite{Pivato-Yassawi-2002I} and \cite{Pivato-Yassawi-2002II}.

However, the uniform Bernoulli measure is not the only $(F,\s)$-invariant measure, indeed every uniform measure
supported on a $(F,\s)$-periodic orbit is $(F,\s)$-invariant. We want to obtain additional conditions which
allow us to characterize the uniform Bernoulli measure. We limit the study to CA which have algebraic and strong
combinatorial properties: the algebraic bipermutative CA. Let $(\az,F)$ be a bipermutative algebraic CA; we
examine the conditions that force an $(F,\s)$-invariant measure $\mu$ to be the Haar measure of $\az$, denoted
by $\lambda_{\az}$. When $\az$ is an infinite product of the finite group $\A$, the Haar measure is the uniform
Bernoulli measure. B. Host, A. Maass and S. Mart\'{\i}nez take this direction in \cite{Host-Maass-Martinez-2003}
and characterize $(F,\s)$-invariant measure of affine bipermutative CA of radius $1$ when the alphabet is $\Z
/p\Z$ with $p$ prime. They show two theorems with different assumptions on the measure $\mu$. M. Pivato gives in
\cite{Pivato-2003} an extension of the first one considering a larger class of algebraic CA but with extra
conditions on the measure and the kernel of $F$. The main result in the present paper provides a generalization
of the second theorem of \cite{Host-Maass-Martinez-2003} which also generalizes Pivato's result.

To introduce more precisely previous work and this article, we need to provide definitions and introduce some
classes of CA. Let $\A$ be a finite set and $\M=\N$ or $\Z$. We consider $\am$, the configuration space of
$\M$-indexed sequences in $\A$. If $\A$ is endowed with the discrete topology, $\am$ is compact and totally
disconnected in the product topology. The {\em shift} map $\s:\am\rightarrow\am$ is defined by $\s(x)_i=x_{i+1}$
for $x=(x_m)_{m\in\M}\in\am$ and $i\in\M$. Denote by $\A^{\ast}$ the set of all finite sequences or words
$w=w_0...w_{n-1}$ with letters in $\A$; by $|w|$ we mean the length of $w\in\A^{\ast}$. Given $w\in\A^{\ast}$
and $i\in\M$, the {\em cylinder set} starting at coordinate $i$ with the word $w$ is $[w]_i=\{x\in\am  :
x_{i,i+|w|-1}=w\}$, the cylinder set starting at $0$ is simply denoted by $[w]$.

A {\em cellular automaton} (CA) is a pair $(\am,F)$ where $\am$ is called the {\em configuration space}, and
$F:\am\to\am$ is a continuous function which commutes with the shift. We can therefore consider $(F,\s)$ as a
$\N\times\M$-action. By Hedlund's theorem \cite{Hedlund-1969}, it is equivalent to give a local function which
acts uniformly and synchronously on the configuration space, that is to say, there is a finite segment
$\U\subset\M$ (named {\em neighborhood}) and a {\em local rule} $\overline{F}:\A^{\U}\rightarrow\A$, such that
$F(x)_m=\overline{F}((x_{m+u})_{u\in\U})$ for all $x\in\am$ and $m\in\M$. The {\em radius} of $F$ is
$r(F)=\max\{|u|: u\in\U\}$; when $\U$ is as small as possible, it is called the {\em smallest neighborhood}. If
the smallest neighborhood is reduced to one point we say that $F$ is {\em trivial}.

Let $\Ba$ be the Borel sigma-algebra of $\am$, we denote by $\Me(\am)$ the set of probability measures on $\am$
defined on the sigma-algebra $\Ba$. As usual, $\s\mu$ (respectively $F\mu$) denotes the measure given by
$\s\mu(B)=\mu(\s^{-1}(B))$ (respectively $F\mu(B)=\mu(F^{-1}(B))$) for $B$ a Borel set; this allows us to
consider the $(F,\s)$-action on $\Me(\am)$. We say that $\mu\in\Me(\am)$ is {\em $\s$-invariant} (respectively
{\em $F$-invariant}) iff $\s\mu=\mu$ (respectively $F\mu=\mu$); obviously $\mu$ is {\em $(F,\s)$-invariant} iff
$\mu$ is $\s$-invariant and $F$-invariant. We denote $\I_{\mu}(\s)=\{B\in\Ba : \mu ( \s^{-1}(B)\triangle B)=
0\}$ the algebra of $\s$-invariant sets $\textrm{mod } \mu$. If $\am$ has a group structure and $\gs$ is a
closed $\s$-invariant subgroup of $\am$, the {\em Haar measure} on $\gs$, denoted $\lambda_{\gs}$, is the unique
measure in $\Me(\am)$ with $\supp(\mu)\subset\gs$ which is invariant by the action of $\gs$. We can characterize
$\lambda_{\gs}$ using {\em characters} in $\widehat{\am}$, which are continuous morphisms from $\am$ to $\C$:
indeed, $\mu=\lambda_{\gs}$ iff $\supp(\mu)\subset \gs$ and $\mu(\chi)=0$ for all $\chi\in\widehat{\am}$ such
that $\chi(\gs)\ne\{1\}$, see~\cite{Guichardet-1968} for more detail. If $\A$ is a finite group and $\am$ is a
product group, the Haar measure of $\am$ corresponds to the {\em uniform Bernoulli measure} defined on a
cylinder set $[u]_i$ by:
$$\lambda_{\am}([u]_i)=\frac{1}{|\A|^{|u|}}.$$

Let $(\am,F)$ be a CA of smallest neighborhood $\U=[r,s]=\{r,...,s\}$. $F$ is {\em left-permutative} iff for any
$u\in\A^{s-r}$ and $b\in\A$, there is a unique $a\in\A$ such that $\overline{F}(au)=b$; $F$ is {\em
right-permutative} iff for any $u\in\A^{s-r}$ and $b\in\A$ there is a unique $a\in\A$ such that
$\overline{F}(ua)=b$. $F$ is {\em bipermutative} iff it is both left and right permutative.

If $\am$ has a topological group structure and if $\s:\am\rightarrow\am$ is a continuous group endomorphism,
$\am$ is called a {\em group shift}. By Hedlund's Theorem~\cite{Hedlund-1969}, the $\s$-commuting multiplication
operator is given by a local rule $\overline{\ast}:\A^{[r,s]}\times\A^{[r,s]}\to\A$. We refer
to~\cite{Kitchens-1987} for more details. If $\am$ is an Abelian group shift and $F:\am\rightarrow\am$ is a
group endomorphisms which commutes with $\s$, then the CA $(\am,F)$ is called {\em algebraic}. If $\A$ has an
Abelian group structure, $\am$ is a compact Abelian group. We say that $(\am,F)$ is a {\em linear} CA if $F$ is
a group endomorphism  or equivalently if $\overline{F}$ is a morphism from $\A^{\U}$ to $\A$. In this case $F$
can be written:
$$ F=\sum_{u\in\U} f_u \circ\s^u \qquad  $$
where for all $u\in\U$, $f_u$ is an endomorphism of $\A$ which is extended coordinate by coordinate to $\am$. We
can write $F$ as a polynomial of $\s$, $F=P_F(\s)$, where $P_F \in\Hom(\A)[X,X^{-1}]$. If $\A=\Z/n\Z$, then an
endomorphism of $\A$ is the multiplication by an element of $\Z/n\Z$. We say that $(\am,F)$ is an {\em affine}
CA if there exists $(\am,G)$ a linear CA and a constant $c\in\am$ such that $F=G+c$. The constant must be
$\s$-invariant.

A linear CA $(\am,F)$ where $F=\sum_{u\in[r,s]}f_u\circ\s^u$ is left (right) permutative of smallest
neighborhood $[r,s]$ if $f_r$ ($f_s$) is a group automorphism. An affine CA $(\am,F+c)$, where $(\am,F)$ is
linear and $c\in\am$, is bipermutative if $(\am,F)$ is bipermutative. So if $\A=\Z/p\Z$ where $p$ is prime, then
any nontrivial affine CA is bipermutative. However, if $p$ is composite, then $F$ is left (right) permutative
iff the leftmost (rightmost) coefficient of $\overline{F}$ is relatively prime to $p$.

Now we can recall the first theorem of \cite{Host-Maass-Martinez-2003}:

\begin{theorem}[\cite{Host-Maass-Martinez-2003}]\label{HMM1}
Let $(\az,F)$ be an affine bipermutative CA of smallest neighborhood $\U=[0,1]$ with $\A=\Z/p\Z$, where $p$ is
prime, and let $\mu$ be an $(F,\s)$-invariant probability measure. Assume that:
\begin{enumerate}
\item $\mu$ is ergodic for $\s$;

\item the measure entropy of $F$ is positive ($h_{\mu}(F)>0$).
\end{enumerate}
Then $\mu=\lambda_{\az}$.
\end{theorem}

The second theorem of \cite{Host-Maass-Martinez-2003} relaxes the $\s$-ergodicity into $(F,\s)$-ergodicity
provided the measure satisfies a technical condition on the sigma-algebra of invariant sets for powers of $\s$:

\begin{theorem}[\cite{Host-Maass-Martinez-2003}]\label{HMM2}
Let $(\az,F)$ be an affine bipermutative CA of smallest neighborhood $\U=[0,1]$ with $\A=\Z/p\Z$ where $p$ is
prime, and let $\mu$ be an $(F,\s)$-invariant probability measure. Assume that:
\begin{enumerate}
\item $\mu$ is ergodic for the $\N\times\Z$-action $(F,\s)$;

\item $\I_{\mu}(\s)=\I_{\mu}(\s^{p(p-1)})$  $\mod\mu$;

\item $h_{\mu}(F)>0$.
\end{enumerate}
Then $\mu=\lambda_{\az}$.
\end{theorem}

M. Pivato gives in \cite{Pivato-2003} a result similar to Theorem \ref{HMM1}, which applies to a larger class of
algebraic CA but with extra conditions on the measure and $\Ker(F)$:

\begin{theorem}[\cite{Pivato-2003}]\label{P2}
Let $\az$ be any Abelian group shift, let $(\az,F)$ be an algebraic bipermutative CA of smallest neighborhood
$\U=[0,1]$ and let $\mu$ be an $(F,\s)$-invariant probability measure. Assume that:
\begin{enumerate}
\item $\mu$ is totally ergodic for $\s$;

\item $h_{\mu}(F)>0$;

\item $\Ker(F)$ contains no nontrivial $\s$-invariant subgroups.
\end{enumerate}
Then $\mu=\lambda_{\az}$.
\end{theorem}

It is possible to extend Theorem~\ref{P2} to a nontrivial algebraic bipermutative CA without restriction on the
neighborhood. In Section~2 of this paper we give entropy formulas for bipermutative CA without restrictions on
the neighborhood. These formulas are the first step to adapt the proof of Theorem~\ref{HMM2} in Section~3 in
order to obtain our main result:

\vspace{6pt} \noindent {\bf Theorem \ref{G1}. } {\em Let $\az$ be any Abelian group shift, let $(\az,F)$ be a
nontrivial algebraic bipermutative CA, let $\gs$ be a closed $(F,\s)$-invariant subgroup of $\az$, let $k\in\N$
such that every prime factor of $|\A|$ divides $k$ and let $\mu$ be an $(F,\s)$-invariant probability measure on
$\az$ with $\supp(\mu)\subset\gs$. Assume that:
\begin{enumerate}
\item $\mu$ is ergodic for the $\N\times\Z$-action $(F,\s)$;

\item $\I_{\mu}(\s)=\I_{\mu}(\s^{kp_1})$ with $p_1$ the smallest common period of all elements of $\Ker(F)$;

\item $h_{\mu}(F)>0$;

\item every $\s$-invariant infinite subgroup of $\ds_{\infty}(F)=\cup_{n\in\N} \Ker(F^n)\cap\gs$ is dense in
$\gs$.
\end{enumerate}
Then $\mu=\lambda_{\gs}$.}

\vspace{5pt} \noindent Theorem~\ref{G1} is a common generalization of Theorem~\ref{HMM2} and Theorem~\ref{P2}
when $\A$ is a cyclic group and $\az$ is the product group. To obtain a generalization of Theorem~\ref{P2} for
any Abelian group $\az$, we must take a weaker assumption for $\ds_{\infty}$, however we need a further
restriction for the probability measure:

\vspace{6pt} \noindent {\bf Theorem \ref{G2}. } {\em Let $\az$ be any Abelian group shift, let $(\az,F)$ be a
nontrivial algebraic bipermutative CA, let $\gs$ be a closed $(F,\s)$-invariant subgroup of $\az$, let $k\in\N$
such that every prime factor of $|\A|$ divides $k$ and let $\mu$ be an $(F,\s)$-invariant probability measure on
$\az$ with $\supp(\mu)\subset\gs$. Assume that:
\begin{enumerate}
\item $\mu$ is ergodic for $\s$;

\item $\I_{\mu}(\s)=\I_{\mu}(\s^{kp_1})$ with $p_1$ the smallest common period of all elements of $\Ker(F)$;

\item $h_{\mu}(F)>0$;

\item every $(F,\s)$-invariant infinite subgroup of $\ds_{\infty}(F)=\cup_{n\in\N} \Ker(F^n)\cap\gs$ is dense in
$\gs$.
\end{enumerate}
Then $\mu=\lambda_{\gs}$.}

\vspace{5pt} \noindent To do this some technical work is required on each of the assumptions. Presently we do
not know how to obtain a common generalization of Theorems~\ref{G1} and~\ref{G2}.

In Section 4 we show how to replace and relax some assumptions of Theorems~\ref{G1} and~\ref{G2}, in particular
how one obtains Theorems~\ref{HMM2} and~\ref{P2} as consequences. First we replace the assumption of positive
entropy of $F$ by the positive entropy of $F^n\circ\s^m$ for some $(n,m)\in\N\times\Z$. Then we give a necessary
and sufficient condition for $D^{\gs}_{\infty}$ to contain no nontrivial $(F,\s)$-invariant infinite subgroups.
This condition is implied by the assumption that $\Ker(F)$ contains no nontrivial $\s$-invariant subgroups.

In Section 5 we restrict the study to linear CA and obtain rigidity results which cannot be deduced from
Theorem~\ref{HMM2} and~\ref{P2}. For example, in Subsection 5.1, we can see that Theorem~\ref{G1} works for
$F=P_F(\s)$ any nontrivial linear CA on $(\Z/p\Z)^{\Z}$ with $p$ prime. In this case Theorem~\ref{HMM2} works
only for CA of radius 1 and Pivato's result works only if $P_F$ is irreducible on $\Z/p\Z$. In Section 6 we give
an application of this work. We stray from the algebraic bipermutative CA case and show measure rigidity for
some affine one-sided invertible expansive CA (not necessary bipermutative) with the help of previous results.

\section{Entropy formulas for bipermutative CA}

Let $(\az,F)$ be a CA, $\Ba$ be the Borel sigma-algebra of $\az$ and $\mu\in\Me(\az)$. We put
$\Ba_n=F^{-n}(\Ba)$ for $n\in\N$. For $\p$ a finite partition of $\az$ and for $\Ba'$ a sub sigma-algebra of
$\Ba$ we denote $H_{\mu}(\p)=-\sum_{A\in\p} \mu(A) \log(\mu(A))$ the entropy of $\p$ and
$H_{\mu}(\p|\Ba')=-\sum_{A\in\p}\int_A\log(\E(1_A|\Ba'))d\mu$ the conditional entropy of $\p$ given $\Ba'$.
Furthermore $h_{\mu}(F)$ denotes the entropy of the measure-preserving dynamical system $(\am,\Ba,\mu,F)$. We
refer to \cite{Petersen-1989} or \cite{Walters-1982} for the definition and main properties.

We define the cylinder partitions $\p=\{[a]:a\in\A\}$ and $\p_{[r,s]}=\{[u]_r : u\in\A^{s-r}\}$. The following
lemma is a more general version of the entropy formula in Lemma $4.3.$ of \cite{Host-Maass-Martinez-2003} (where
this Lemma is proved for CA with radius 1):

\begin{lemma}\label{cG}
Let $(\az,F)$ be a bipermutative CA of smallest neighborhood $\U=[r,s]$ with $r\leq 0 \leq s$ and let $\mu$ be
an $F$-invariant probability measure on $\az$. Then $h_{\mu}(F)=H_{\mu}(\p_{[0,s-r-1]}|\Ba_1)$.
\end{lemma}
\begin{proof}
We have $h_{\mu}(F)=\lim_{l\to\infty}h_{\mu}(F,\p_{[-l,l]})$ with:
$$ h_{\mu}(F,\p_{[-l,l]}) = \lim_{T\to\infty}H_{\mu}(\p_{[-l,l]}|\bigvee_{n=1}^{T}F^{-n}(\p_{[-l,l]})) = H_{\mu}(\p_{[-l,l]}|\bigvee_{n=1}^{\infty}F^{-n}(\p_{[-l,l]})).$$

Let $l\geq s-r$. By bipermutativity of $F$, for $T\geq 1$, it is equivalent to know $(F^n(x)_{[-l,l]})_{n\in
[1,T]}$ and to know $F(x)_{[Tr-l,Ts+l]}$. This means that
$\bigvee_{n=1}^{T}F^{-n}(\p_{[-l,l]})=F^{-1}(\p_{[Tr-l,Ts+l]})$. By taking the limit as $l\to\infty$, we deduce
(with the convention $\infty.0=0$):
$$ \bigvee_{n=1}^{\infty}F^{-n}(\p_{[-l,l]})=F^{-1}(\p_{[\infty.r-l,\infty.s+l]}).$$
So we have:
$$ h_{\mu}(F,\p_{[-l,l]}) =  H_{\mu}(\p_{[-l,l]}|F^{-1}\p_{[\infty.r-l,\infty.s+l]}) .$$

Similarly, by bipermutativity of $F$, the knowledge of $F(x)_{[\infty.r-l,\infty.s+l]}$ and $x_{[0,s-r-1]}$
allows us to know $x_{[-l,l]}$ and vice versa. We deduce:
$$ \p_{[0,s-r-1]}\vee F^{-1}(\p_{[\infty.r-l,\infty.s+l]}) = \p_{[-l,l]}\vee F^{-1}(\p_{[\infty.r-l,\infty.s+l]}) .$$
Therefore,
$$  h_{\mu}(F,\p_{[-l,l]}) =  H_{\mu}(\p_{[0,s-r+1]}|F^{-1}(\p_{[\infty.r-l,\infty.s+l]})) .$$

If $r<0<s$, then $\p_{[\infty.r-l,\infty.s+l]}=\Ba_1$. Otherwise, by taking the limit as $l\to\infty$ and using
the martingale convergence theorem, we obtain $h_{\mu}(F)=H_{\mu}(\p_{[0,s-r-1]}|\Ba_1)$.
\end{proof}

When $\mu$ is an $(F,\s)$-invariant probability measure, it is possible to express the entropy of a
right-permutative CA according to the entropy of $\s$.

\begin{prop}\label{entropyperm}
Let $(\az,F)$ be a right-permutative CA of neighborhood $\U=[0,s]$, where $s$ is the smallest possible value and
let $\mu$ be an $(F,\s)$-invariant probability measure. Then $h_{\mu}(F)=s\, h_{\mu}(\s)$.
\end{prop}
\begin{proof}
Let $N\in\N$ and $l\geq s$. By right-permutativity, since $\U=[0,s]$, for all $x\in\az$ it is equivalent to know
$(F^n(x)_{[-l,l]})_{n\in [0,N]}$ and $x_{[-l,l+Ns]}$; this means that: $$\bigvee_{n=0}^N F^{-n}(\p_{[-l,l]}) =
\p_{[-l,l+Ns]}.$$ So for $l\geq s$ we have:
\begin{eqnarray*}
h_{\mu}(F,\p_{[-l,l]}) & = &  \lim_{N\to\infty} \frac{1}{N} H_{\mu}(\bigvee_{n=0}^N F^{-n}(\p_{[-l,l]})) \\
& = & \lim_{N\to\infty}\frac{1}{N} H_{\mu}(\p_{[-l,l+Ns]})  \\
&=& \lim_{N\to\infty} - \frac{1}{N} \sum_{u\in \A^{Ns+2l}} \mu([u])\log(\mu[u]) \\
& =& \lim_{N\to\infty} -\frac{Ns+2l}{N} \frac{1}{Ns+2l} \sum_{u\in \A^{Ns+2l}} \mu([u])\log(\mu[u])\\
& =& s \, h_{\mu}(\s).
\end{eqnarray*}
We deduce that $h_{\mu}(F)=\lim_{l\to\infty } h_{\mu}(F,\p_{[-l,l]}) = s\, h_{\mu}(\s)$.
\end{proof}
\begin{remark}
We have a similar formula for a left-permutative CA of neighborhood $\U=[r,0]$. Moreover, it is easy to see that
this proof is true for a right-permutative CA on $\an$.
\end{remark}

\begin{cor}\label{centropie}
Let $(\az,F)$ be a bipermutative CA of smallest neighborhood $\U=[r,s]$, and let $\mu$ be an $(F,\s)$-invariant
probability measure on $\az$. We have:
$$h_{\mu}(F)=
\begin{cases}
s\, h_{\mu}(\s) & \textrm{ if  }  s\geq r \geq 0,\\
(s-r)\, h_{\mu}(\s) & \textrm{ if  }  s\geq 0\geq r,\\
-r \, h_{\mu}(\s) & \textrm{ if  } 0\geq s \geq r.
\end{cases}$$
\end{cor}
\begin{proof}
Cases where $s\geq r \geq 0$ or $0\geq s \geq r$ can be directly deduced from Proposition~\ref{entropyperm}.

When $s\geq 0\geq r$, the CA $(\az,\s^{-r}\circ F)$ is bipermutative of smallest neighborhood $[0,s-r]$. Since
$\s$ is bijective, we deduce that $\Ba$ is $\s$-invariant. Thus, $F^{-1}(\Ba)=(\s^{-r}\circ F)^{-1}(\Ba)$. Since
$\mu$ is $(F,\s)$-invariant, by Lemma~\ref{cG}, one has:
$$h_{\mu}(F)=H_{\mu}(\p_{[0,s-r-1]}|F^{-1}(\Ba))=H_{\mu}(\p_{[0,s-r-1]}|(\s^{-r}\circ
F)^{-1}(\Ba))=h_{\mu}(\s^{-r}\circ F).$$

The result follows from Proposition~\ref{entropyperm}.
\end{proof}

\begin{remark}
It is not necessary to use Lemma~\ref{cG}. Corollary~\ref{centropie} can be proved by a similar method of
Proposition~\ref{entropyperm}.
\end{remark}

A bipermutative CA $(\az,F)$ of smallest neighborhood $\U$ is topologically conjugate to $((\A^t)^{\N},\s)$
where $t=\max(\U\cup\{0\})-\min(\U\cup\{0\})$, via the conjugacy $\varphi:x\in\az\to (F(x)_{[0,t]})_{n\in\N}$.
So the uniform Bernoulli measure is a maximal entropy measure. Thus from Corollary \ref{centropie} we deduce an
expression of $h_{\textrm{top}}(F)$. This implies a result of \cite{Ward-2000} which compute the topological
entropy for linear CA on $(\Z /p\Z )^{\Z}$ with $p$ prime by algebraic methods. Moreover this formula gives
Lyapunov exponents for permutative CA according to the definition of \cite{Shereshevsky2-1992} or
\cite{Tisseur-2000}.

\section{Proof of main theorems}

Now we consider $(\az,F)$ a bipermutative algebraic CA of smallest neighborhood $\U=[r,s]$. For $y\in\az$ call
$T_y$ the translation $x\mapsto x+y$ on $\az$. For every $n\in\N$, we write $D_n(F)=\Ker (F^n)$; if there is no
ambiguity we just denote it by $D_n$. Clearly $D_n$ is a subgroup of $D_{n+1}$. Denote $\partial
D_{n+1}=D_{n+1}\setminus D_n$ for all $n\in\N$. By bipermutativity we have $|D_n|=|D_1|^n=|\A|^{(s-r)n}$ where
$|.|$ denotes the cardinality of the set. We can consider the subgroup $D_{\infty}(F)=\cup_{n\in\N}D_n(F)$ of
$\az$, we denote it by $D_{\infty}$ if there is no ambiguity; it is dense in $\az$ since $F$ is bipermutative.
Every $D_n$ is finite and $\s$-invariant so every $x\in D_n$ is $\s$-periodic. Let $p_n$ be the smallest common
period of all elements of $D_n$. Then $p_n$ divides $|D_n|!$.

Let $\Ba$ be the Borel sigma-algebra of $\az$ and let $\mu$ be a probability measure on $\az$. Put
$\Ba_n=F^{-n}(\Ba)$ for every $n\in\N$, it is the sigma-algebra generated by all cosets of $D_n$. For every
$n\in\N$ and $\mu$-almost every $x\in\az$, the conditional measure $\mu_{n,x}$ is defined for every measurable
set $U\subset \az$ by $\mu_{n,x}(U)=\E(\mathbf{1}_U|\Ba_n)(x)$. Its main properties are:
\begin{itemize}

\item[(A)] For $\mu$-almost every $x\in\az$, $\mu_{n,x}$ is a probability measure on $\az$ and
$\supp(\mu_{n,x})\subset F^{-n}(\{F^n(x)\})=x+D_n$.

\item[(B)] For all measurable sets $U\subset\az$, the function $x\rightarrow\mu_{n,x}(U)$ is $\Ba_n$-measurable
and $\mu_{n,x}=\mu_{n,y}$ for every $y\in F^{-n}(\{F^n(x)\})=x+D_n$.

\item[(C)] Let $G:\az\to\az$ be a measurable map and let $U$ be a measurable set. For $\mu$-almost every
$x\in\az$ one has $\E(\mathbf{1}_{G^{-1}(U)}|G^{-1}(\Ba))(x)=\E(\mathbf{1}_U|\Ba)(G(x))$. So
$\s^m\mu_{n,x}=\mu_{n,\s^m(x)}$ and $F \mu_{n+1,x}=\mu_{n,F(x)}$ for $\mu$-almost every $x\in\az$ and every
$n\in\N$.

\item[(D)] Since $\Ba_n$ is $T_d$-invariant for $d\in D_n$, by (C) one has $\mu_{n,x}=\mu_{n,x+d}$.
\end{itemize}

For all $n\in\N$ define $\zeta_{n,x}=T_{-x}\mu_{n,x}$; it is a probability measure concentrated on $D_n$. The
previous four properties of conditional measures can be transposed to $\zeta_{n,x}$:
\begin{lemma}\label{zetaprop}
Fix $n\in\N$. For $\mu$-almost all $x\in\az$, the following are true:
\begin{itemize}
\item[(a)] $\zeta_{n,x+d}=T_{-d}\zeta_{n,x}$ for every $d\in D_n$.

\item[(b)] $\s^m\zeta_{n,x}=\zeta_{n,\s^m(x)}$ for every $m\in\Z$ and $F\zeta_{n+1,x}=\zeta_{n,F(x)}$.

\item[(c)] For every $m\in p_n\Z$, we have $\s^{m}\zeta_{n,x}=\zeta_{n,x}$. Hence $x\to\zeta_{n,x}$ is
$\s^m$-invariant.
\end{itemize}
\end{lemma}
\begin{proof}
(a) is by Property (D). (b) is by Property (C). And (c) is because $\supp(\zeta_{n,x})\subset D_n$.
\end{proof}

For $n>0$ and $d\in D_n$ we define:
$$ E_{n,d}=\{x\in\az : \zeta_{n,x}(\{d\})>0\}\, \textrm{ and }\, E_n=\bigcup_{d\in \partial D_n} E_{n,d}.$$
Then $E_{n,d}$ is $\s^{p_n}$-invariant by Lemma~\ref{zetaprop}(c), and $E_n$ is $\s$-invariant, because
$\partial D_n$ is $\s$-invariant. We write $\eta(x)=\zeta_{1,x}(\{0\})=\mu_{1,x}(\{x\})$. The function $\eta$ is
$\s$-invariant and $E_1=\{x\in\az : \eta(x)<1\}$. Therefore one has:
$$\eta(F^{n-1}(x))=\mu_{1,F^{n-1}(x)}(\{F^{n-1}(x)\})=\mu_{1,F^{n-1}(x)}(F^{n-1}(x+D_{n-1}))\eg{(*)} \mu_{n,x}(x+D_{n-1})=\zeta_{n,x}(D_{n-1}),$$
where $(*)$ is by property (C). Thus $E_n=\{x\in\az : \zeta_{n,x}(D_{n-1})<1\}=F^{-n+1}(E_1)$.

Let $\gs\subset\az$ be a closed $(F,\s)$-invariant subgroup of $\az$. We denote $\ds_n=D_n\cap\gs$ and $\partial
\ds_{n+1}=\ds_{n+1}\setminus \ds_n$ for all $n\in\N$ and $\ds_{\infty}=D_{\infty}\cap\gs$.
\begin{remark}\label{rem}
For $\mu$ an $(F,\s)$-invariant probability measure such that $\supp(\mu)\subset\gs$, we remark that for every
$n\in\N$ and $\mu$-almost every $x\in\az$, $\supp(\mu_{n,x})\subset x+\ds_n\subset\gs$ and
$\supp(\zeta_{n,x})\subset\ds_n$. So for all $n\in\N$ and $d\in \partial D_n$, if $d\notin\gs$ one has
$\mu(E_{n,d})=0$.
\end{remark}

\begin{lemma}\label{lemmainv}
Let $\mu$ be a $\s$-invariant measure on $\az$. If there exist $k\in\N$ such that
$\I_{\mu}(\s)=\I_{\mu}(\s^{k})$ then for all $n\geq 1$ one has  $\I_{\mu}(\s)=\I_{\mu}(\s^{k^n})$.
\end{lemma}
\begin{proof}
Applying the ergodic decomposition theorem to $(\az,\Ba,\mu,\s)$, to prove $\I_{\mu}(\s)=\I_{\mu}(\s^{k^n})$ it
is equivalent to prove that almost every $\s$-ergodic component $\dd$ of $\mu$ is ergodic for $\s^{k^n}$. The
proof is done by induction.

The base case $\I_{\mu}(\s)=\I_{\mu}(\s^{k^1})$ is true by hypothesis. Let $n\geq 2$ and assume that this
property holds for $n-1$ and does not hold for $n$. That is to say we consider a $\s$-ergodic component $\dd$ of
$\mu$ (by induction it is also $\s^{k^{n-1}}$-ergodic) which is not $\s^{k^n}$-ergodic. There exist
$\lambda\in\C$ such that $\lambda^{k^n}= 1 $ and $\lambda^{k^{n-1}}\ne 1$ and a non constant function
$h:\az\to\C$ such that $h(\s(x))=\lambda h(x)$ for $\dd$-almost every $x\in\az$. We deduce that
$h^{k}(\s(x))=\lambda^{k}h^{k}(x)$ and $h^{k}(\s^{k^{n-1}}(x))=\lambda^{k^n}h^{k}(x)=h^{k}(x)$ for $\dd$-almost
every $x\in\az$. By $\s^{k^{n-1}}$-ergodicity of $\dd$, $h^{k}$ is constant $\delta$-almost everywhere, so
$\lambda^{k}=1$ which is a contradiction.
\end{proof}

\begin{remark}\label{annexelemma}
If $k$ divides $k'$ then $\I_{\mu}(\s)\subset\I_{\mu}(\s^k)\subset\I_{\mu}(\s^{k'})$. So if
$\I_{\mu}(\s)=\I_{\mu}(\s^{k'})$ one also has $\I_{\mu}(\s)=\I_{\mu}(\s^k)$.
\end{remark}

We recall the main theorem:
\begin{theorem}\label{G1}
Let $\az$ be any Abelian group shift, let $(\az,F)$ be a nontrivial algebraic bipermutative CA, let $\gs$ be a
closed $(F,\s)$-invariant subgroup of $\az$, let $k\in\N$ such that every prime factor of $|\A|$ divides $k$ and
let $\mu$ be an $(F,\s)$-invariant probability measure on $\az$ with $\supp(\mu)\subset\gs$. Assume that:
\begin{enumerate}
\item $\mu$ is ergodic for the $\N\times\Z$-action $(F,\s)$;

\item $\I_{\mu}(\s)=\I_{\mu}(\s^{k p_1})$ with $p_1$ the smallest common period of all elements of $\Ker(F)$;

\item $h_{\mu}(F)>0$;

\item every $\s$-invariant infinite subgroup of $\ds_{\infty}(F)=\cup_{n\in\N} \Ker(F^n)\cap\gs$ is dense in
$\gs$.
\end{enumerate}
Then $\mu=\lambda_{\gs}$.
\end{theorem}
\begin{proof}
For all $n\in\Z$, $F$ is bipermutative iff $\s^n\circ F$ is bipermutative. Since $F$ is nontrivial, by
Corollary~\ref{centropie}, we deduce that $h_{\mu}(\s^n\circ F)>0$ for all $n\in\Z$. Moreover $\mu$ is
$\s$-invariant. So we can assume that the smallest neighborhood of $F$ is $[0,r]$ with $r\in\N$.

\Claim{For all $n\in\N$, $\I_{\mu}(\s)=\I_{\mu}(\s^{kp_n})$, where $p_n$ is the smallest common $\s$-period of
$D_n$.} \bclaimprf Let $n\in\N$. Every $x\in D_n$ is a $\s$-periodic point of $\s$-period $p_n$, so by
bipermutativity, every $y\in F^{-1}(\{x\})$ is $\s$-periodic. Since $\s^{p_n}(y)\in F^{-1}(\{x\})$, one has that
$p_n$ divides the $\s$-period of $y$. We deduce that $p_n$ divides $p_{n+1}$. Moreover there exists $d\in D_1$
such that $\s^{p_n}(y)=y+d$, so $\s^{|D_1|p_n}(y)=y+|D_1|d=y$. We deduce that $p_{n+1}$ divides $|\A|^r p_n$,
because $|D_1|=|\A|^r$. By induction $p_n$ divides $|\A|^{r(n-1)} p_1$. If $m$ is large enough, then
$|\A|^{r(n-1)}$ divides $k^m$, hence $p_n$ divides $|\A^{r(n-1)}p_1|$ which divides $(kp_1)^m$. Thus,
$\I_{\mu}(\s^{kp_n})=\I_{\mu}(\s)$ by Remark~\ref{annexelemma}, because $\I_{\mu}(\s^{(kp_n)^m})=\I_{\mu}(\s)$
 by Lemma~\ref{lemmainv} and hypothesis (2) of Theorem~\ref{G1}.
\eclaimprf

\Claim{For $n\in\N$ and $d\in D_n$, the measure $T_d(\mathbf{1}_{E_{n,d}}\mu)$ is absolutely continuous with
respect to $\mu$.} \bclaimprf Let $A\in\Ba$ be such that $\mu(A)=0$. Since
$\mu(A)=\int_{\az}\mu_{n,x}(A)d\mu(x)$, we deduce that $\mu_{n,x}(A)=0$ for $\mu$-almost every $x\in\az$. In
particular, $0=\mu_{n,x}(A)\geq \mu_{n,x}(\{x+d\})=\zeta_{n,x}(\{d\})$, for $\mu$-almost every $x\in T_{-d}(A)$
because $x+d\in A$. Thus $x\notin E_{n,d}$ so $\mu(T_{-d}(A)\cap E_{n,d})=0$. This implies that
$T_d(\mathbf{1}_{E_{n,d}}\mu)(A)=0$, so $T_d(\mathbf{1}_{E_{n,d}}\mu)$ is absolutely continuous with respect to
$\mu$. \eclaimprf

To prove the theorem, we consider $\chi\in\widehat{\az}$ with $\mu(\chi)\ne 0$ and we show that $\chi(x)=1$ for
all $x\in\gs$. We consider $\Gamma=\{d\in \ds_{\infty}:\chi(d)=\chi(\s^m(d)),\forall m\in\Z\}$, a $\s$-invariant
subgroup of $\ds_{\infty}$. We want to show that $\Gamma$ is infinite and hence, dense in $\gs$ by hypothesis
(4). From this we will deduce that $\chi$ must be constant.

\Claim{There exists $N\subset\az$ with $\mu(N)=1$ and $F(N)=N$ (up to a set of measure zero), satisfying the
following property: For any $n\in\N$ and $d\in \partial\ds_n$, if there exists $x\in E_{n,d}\cap
 N$ with $\zeta_{n,x}(\chi)\ne 0$, then $d\in\Gamma$.}\bclaimprf For $n\in\N$, the function $x\rightarrow
\zeta_{n,x}$ is $\s^{kp_n}$-invariant by Lemma~\ref{zetaprop}{(c)}. Since $\I_{\mu}(\s)=\I_{\mu}(\s^{kp_n})$ by
Claim 1, we deduce that $\zeta_{n,x}$ is $\s$-invariant. So for $\mu$-almost every $x\in\az$ and for any
$m\in\Z$, we have $\s^m\zeta_{n,x}=\zeta_{n,x} \ (\dag)$. Since $T_d(\mathbf{1}_{E_{n,d}}\mu)$ is absolutely
continuous with respect to $\mu$ by Claim 2, we have $\s^m\zeta_{n,x+d}=\zeta_{n,x+d} \ (\ddag)$ too, for
$\mu$-almost every $x\in E_{n,d}$, for every $d\in D_n$ and for every $m\in\Z$. We can compute:
$$ T_{-\s^m d}\zeta_{n,x} \eg{(\dag)} T_{-\s^m d}\s^m\zeta_{n,x} = \s^m T_{-d}\zeta_{n,x} \eg{(\ast)} \s^m\zeta_{n,x+d} \eg{(\ddag)} \zeta_{n,x+d} \eg{(\ast)} T_{-d}\zeta_{n,x},$$
where $(\dag)$ and $(\ddag)$ are as above, and $(\ast)$ is by Lemma~\ref{zetaprop}(a). So $T_{\s^m
d-d}\zeta_{n,x}=\zeta_{n,x}$ and by integration $(1-\chi(\s^m d-d))\zeta_{n,x}(\chi)=0$ for $\mu$-almost every
$x\in E_{n,d}$. Thus, there exists $N\subset\az$ with $\mu(N)=1$, such that for all $d\in D_n$ and $x\in
E_{n,d}\cap N$, if $\zeta_{n,x}(\chi)\ne 0$, then $\chi(\s^m(d))\chi(d)^{-1}=\chi(\s^m(d)-d)=1$. Hence
$\chi(\s^m(d))=\chi(d)$ for all $m\in\Z$, which means $d\in\Gamma$. Moreover the set $N$ is $F$-invariant up to
a set of measure $0$, because $\mu$ is $F$-invariant, thus $\mu(F(N))=F\mu(F(N))=\mu(F^{-1}(F(N)))\geq\mu(N)=1$.
\eclaimprf

\Claim{There exists $n_0\in\N$ such that, if we define $B=\{x\in N : \E(\chi|\Ba_n)(x)\ne 0 , \forall n\geq
n_0\}$, then $\mu(B)>0$. Moreover for all $n\geq n_0$, and any $d\in\partial\ds_n$, if $E_{n,d}\cap
B\ne\emptyset$, then $d\in\Gamma$.} \bclaimprf One has
$\lim_{n\to\infty}\E(\chi|\Ba_n)=\E(\chi|\cap_{m>1}\Ba_m)$ by the Martingale Convergence Theorem, and this
function is not identically $0$ because its integral is equal to $\mu(\chi)\ne 0$. Thus we can choose $n_0$ such
that $B=\{x\in N : \E(\chi|\Ba_n)(x)\ne 0 , \forall n\geq n_0\}$ satisfies $\mu(B)>0$. Moreover, we have:
$$\E(\chi|\Ba_n)(x)=\int_{\az}\chi d\mu_{n,x} = \chi(x)\zeta_{n,x}(\chi).$$
By Claim 3, for any $n\geq n_0$ and any $d\in \partial\ds_n$, if there is $x\in E_{n,d}\cap B$ then
$d\in\Gamma$. \eclaimprf

\Claim{$\mu(E_1)>0$.} \bclaimprf
 Let $A\in\p_{[0,r-1]}$. Let $x\in A$ and $d\in D_1$ such that $x+d\in A$. One has $x_{[0,r-1]}=(x+d)_{[0,r-1]}$
and $F(x)=F(x+d)$. By bipermutativity, one deduces that $x=x+d$, that is to say $d=0$. Therefore, for any $x\in
A$ and for any $d\in \partial D_1$, we have $x+d \notin A$. Thus, $A\cap F^{-1}(\{F(x)\})=A\cap (x+D_1)=\{x\}$.
Thus, (A) implies that $\E(\mathbf{1}_A|\Ba_1)(x)=\mu_{1,x}(A)= \mu_{1,x}(\{x\})=\eta(x)$. By Lemma \ref{cG},
\begin{eqnarray*}
h_{\mu}(F) & = & H_{\mu}(\p_{[0,r-1]}|\Ba_1) \\
&=& - \sum_{A\in \p_{[0,r-1]}} \int_A \log\big( \E(\mathbf{1}_A|\Ba_1)\big)d\mu\\
&=& \int_{\az} -\log(\eta(x))d\mu(x) \\
&\leeeq{(\ast)}& \int_{E_1}-\log(\eta(x))d\mu(x),
\end{eqnarray*}
where $(\ast)$ is because $E_1=\{x\in\az:\eta(x)<1\}$. But $h_{\mu}(F)>0$ by hypothesis $(3)$. This proves
Claim~5. \eclaimprf

\Claim{$\Gamma$ is infinite.} \bclaimprf For $\mu$-almost every $x\in\az$ one has:
$$  \frac{1}{n}\sum_{j=1}^{n+1}\mathbf{1}_{E_j}(x)\eg{(1)}\frac{1}{n}\sum_{j=1}^{n+1}\mathbf{1}_{F^{-j+1}(E_1)}(x) = \frac{1}{n}\sum_{j=0}^{n}\mathbf{1}_{E_1}(F^j(x))\eg{(2)} \frac{1}{n^2}\sum_{j,k=0}^{n}\mathbf{1}_{E_1}(\s^kF^j(x))\rightarrrow{(3)} \mu(E_1)\geeeq{(4)}0 .$$
Here, (1) is because $E_j=F^{-j+1}(E_1)$ for all $j\in\N$, (2) is because $E_1$ is $\s$-invariant, (3) is the
Ergodic Theorem and hypothesis~$1$, and~(4) is by Claim~5.

It follows that for $\mu$-almost every $x\in\az$, there are infinitely many values of $n>0$ such that $x\in
E_n$. Thus $\mu(\bigcap_{m\in\N}\bigcup_{n\geq m}E_{n})=1$. Since $\mu(B)>0$, we deduce that
$\mu(\bigcap_{m\in\N}\bigcup_{n\geq m}E_{n}\cap B)>0$. For all $n\in\N$, if $d\notin\supp(\mu)\subset\gs$, then
Remark~\ref{rem} implies that $\mu(E_{n,d})=0$. We can conclude that $\{d\in \ds_{\infty}:\exists n\in\N
\textrm{ such that }d\in
\partial D_n \textrm{ and } E_{n,d}\cap B\ne\emptyset\}$ is infinite and by Claim 4, it is a subset of $\Gamma$.
Therefore $\Gamma$ is infinite. \eclaimprf

If we consider $\Gamma'=(\id_{\az}-\s)\Gamma$, we have an infinite $\s$-invariant subgroup of $\ds_{\infty}$
because $\Ker(\id_{\az}-\s)$ is finite. Hypothesis~(4) then implies that $\Gamma'$ is dense in $\gs$, but by
construction, $\chi(\Gamma')=\{1\}$, so by continuity of $\chi$, $\chi(x)=1$ for all $x\in\gs$.
Contrapositively, we must have $\mu(\chi)=0$ for all $\chi\in\widehat{\az}$ such that $\chi(\gs)\ne\{1\}$. Since
$\supp(\mu)\subset\gs$, we conclude that $\mu=\lambda_{\gs}$.
\end{proof}

\begin{remark}
The proof of this theorem works if $(\an,F)$ is a right-permutative algebraic CA where all $x\in D_1=\Ker(F)$
 are $\s$-periodic, but this last assumption is possible only if $F$ is also left-permutative, therefore it is a false generalization.
\end{remark}

\begin{remark}
Let $(\az,F)$ be a nontrivial algebraic bipermutative CA and let $\gs$ be a closed $(F,\s)$-invariant subgroup
of $\az$ which verifies hypothesis~(4) of Theorem~\ref{G1}. Let $c\in\gs$ be a $\s$-invariant configuration. We
define the CA $G=F+c$. Let $\mu$ be a $(G,\s)$-invariant probability measure on $\az$. If $\mu$ verifies the
assumptions of Theorem~\ref{G1} for the $\N\times\Z$-action induced by $(G,\s)$, then $\mu=\lambda_{\gs}$.
\end{remark}

Assumption (4) becomes more natural when it is replaced by ``every $(F,\s)$-invariant infinite subgroup of
$\ds_{\infty}(F)=\cup_{n\in\N} \Ker(F^n)\cap\gs$ is dense in $\gs$''. It is not clear that this condition is
implied by the assumptions of Theorem~\ref{G1}. However if we consider a $\s$-ergodic measure we can prove:

\begin{theorem}\label{G2}
Let $\az$ be any Abelian group shift, let $(\az,F)$ be a nontrivial algebraic bipermutative CA, let $\gs$ be a
closed $(F,\s)$-invariant subgroup of $\az$, let $k\in\N$ such that every prime factor of $|\A|$ divides $k$ and
let $\mu$ be an $(F,\s)$-invariant probability measure on $\az$ with $\supp(\mu)\subset\gs$. Assume that:
\begin{enumerate}
\item $\mu$ is ergodic for $\s$;

\item $\I_{\mu}(\s)=\I_{\mu}(\s^{k p_1})$ with $p_1$ the smallest common period of all elements of $\Ker(F)$;

\item $h_{\mu}(F)>0$;

\item every $(F,\s)$-invariant infinite subgroup of $\ds_{\infty}(F)=\cup_{n\in\N} \Ker(F^n)\cap\gs$ is dense in
$\gs$.
\end{enumerate}
Then $\mu=\lambda_{\gs}$.
\end{theorem}
\begin{proof}
A measure $\s$-ergodic is $(F,\s)$-ergodic so results from Claim~1 to Claim~6 hold.

\Claim{Let $B'=\cup_{j\in\Z}\s^j(\{x\in N : \E(\chi|\Ba_n)(x)\ne 0, \forall n\in\N \})$. Then $\mu(B')=1$.}
\bclaimprf By Claim~4, $\mu(B)>0$ where $B=\{x\in N:\E(\chi|\Ba_n)(x)\ne 0, \forall n\geq n_0\}$. Thus there
exists $k\in[0,3]$ such that $B_{n_0}=\{x\in N:\Re(\textbf{i}^k\E(\chi|\Ba_{n_0}))>0, \forall n\geq n_0\}$
verifies $\mu(B_{n_0})>0$, where $\textbf{i}^2=-1$. Since $B_{n_0}\in\Ba_{n_0}\subset\Ba_{n_0-1}$, one has:
$$\int_{B_{n_0}}\Re(\textbf{i}^k\E(\chi|\Ba_{n_0-1}))(x)d\mu=\int_{B_{n_0}}\Re(\textbf{i}^k\E(\chi|\Ba_{n_0})(x))d\mu>0$$
So $B_{n_0-1}=\{x\in B_{n_0}:\Re(\textbf{i}^k\E(\chi|\Ba_{n_0-1})(x))>0\} = \{x\in N:
\Re(\textbf{i}^k\E(\chi|\Ba_{n})(x))>0, \forall n\geq n_0-1\}$ verify $\mu(B_{n_0-1})>0$. By induction
$\mu(B_0)>0$, so $\mu(B')>0$. Since $B'$ is $\s$-invariant, $\mu(B')=1$ by $\s$-ergodicity from hypothesis~(1).
\eclaimprf

\Claim{Let $n\in\N$ and let $d\in\partial\ds_n$. If $E_{n,d}\cap B'$ is nonempty then
$d\in\Gamma=\{d'\in\ds_{\infty}:\chi(d')=\chi(\s^m(d')), \forall m\in\Z\}$} \bclaimprf Let $d\in\partial\ds_n$
and let $x\in E_{n,d}\cap B'$. There exists $j\in\Z$ such that:
$$0\ne \E(\chi|\Ba_n)(\s^j(x))=\int_{\az}\chi d\mu_{n,\s^j(x)}=\chi(\s^j(x))\zeta_{n,\s^j(x)}(\chi)\eg{(\ast)}\chi(\s^j(x))\zeta_{n,x}(\chi).$$
Here $(\ast)$ is because $x\to\zeta_{n,x}$ is $\s^{kp_n}$-invariant by Lemma~\ref{zetaprop}{(c)} and
$\I_{\mu}(\s)=\I_{\mu}(\s^{kp_n})$ by Claim~1, so $x\to\zeta_{n,x}$ is $\s$-invariant. One deduces that
$\zeta_{n,x}(\chi)\ne 0$. But $x\in E_{n,d}\cap N$, so $d\in\Gamma$ by Claim~3. \eclaimprf

\Claim{Let $n\geq 1$ and let $d\in\partial\ds_n$. For $\mu$-almost all $x\in E_{n,d}\cap B'$ one has $F(x)\in
E_{n-1,F(d)}\cap B'$.} \bclaimprf Let $d\in\partial\ds_n$ and $x\in E_{n,d}\cap B'$. One has:
$$\zeta_{n-1,F(x)}(\{F(d)\}) \eg{(1)} \zeta_{n,x}(F^{-1}(\{F(d)\}))\geq \zeta_{n,x}(\{d\})\geeeq{(2)}0.$$
Here (1) is by Lemma~\ref{zetaprop}{(b)} and (2) is because $x\in E_{n,d}$. We deduce that $F(x)\in
E_{n-1,F(d)}$. Since $\mu(B')=1$ by Claim~7 and $\mu$ is $F$-invariant, one has $\mu(\cap_{n\in\N}F^{-n}(B'))=1$
so $F(x)\in E_{n-1,F(d)}\cap B'$ for $\mu$-almost all $x\in E_{n,d}\cap B'$.  \eclaimprf

\Claim{$\cap_{n\in\N} F^{-n}\Gamma$ is infinite.} \bclaimprf Let $n\geq 0$. The set $E_n=F^{-n+1}(E_1)$ is
$\s$-invariant since $E_1$ is $\s$-invariant and $F$ commutes with $\s$. Moreover $\mu(E_n)=\mu(E_1)>0$ by
Claim~5. By $\s$-ergodicity (hypothesis $(1)$), $\mu(E_n)=1$ so $\mu(E_n\cap B')=1$ by Claim~7. For all $n\geq
1$, there exists $d_n\in\partial\ds_n$ such that $\mu(E_{n,d_n}\cap B')>0$, and thus, by Claim~9,
$\mu(E_{n-k,F^k(d_n)}\cap B')>0$ for all $k\in [0,n]$. That is to say $F^k(d_n)\in\Gamma$ for $k\in[0,n]$ by
Claim~8. One deduces that $\cap_{n\in\N} F^{-n}\Gamma$ is infinite since it contains $d_n$ for all
$n\in\N$.\eclaimprf

If we consider $\Gamma''=(\id_{\az}-\s)(\cap_{n\in\N}F^{-n}\Gamma)$, we have an infinite $(F,\s)$-invariant
subgroup of $\ds_{\infty}$ because $\Ker(\id_{\az}-\s)$ is finite. We deduce that $\Gamma''$ is dense in $\gs$
by condition (4), but $\chi(\Gamma'')=\{1\}$ by construction, so by continuity of $\chi$, $\chi(x)=1$ for all
$x\in\gs$. Contrapositively, we must have $\mu(\chi)=0$ for all $\chi\in\widehat{\az}$ such that
$\chi(\gs)\ne\{1\}$. Since $\supp(\mu)\subset\gs$, we conclude that $\mu=\lambda_{\gs}$.
\end{proof}

\begin{cor}\label{cfactor}
Let $\az$ be any Abelian group shift, let $(\az,F)$ be an algebraic bipermutative CA. Let $\gs$ be a closed
$(F,\s)$-invariant subgroups of $\az$ such that there exists $\pi:\az\to\gs$ a surjective continuous morphism
which commutes with $F$ and $\s$ ($(\gs,\s,F)$ is a dynamical and algebraic factor of $(\az,\s,F)$). Let
$k\in\N$ be such that every prime factor of $|\A|$ divides $k$. Let $\mu$ be an $(F,\s)$-invariant probability
measure on $\az$. Assume that:
\begin{enumerate}
\item $\mu$ is ergodic for the $\N\times\Z$-action $(F,\s)$;

\item $\I_{\mu}(\s)=\I_{\mu}(\s^{kp_1})$ with $p_1$ the smallest common period of all elements of $\Ker(F)$;

\item $h_{\pi\mu}(F)>0$;

\item every $\s$-invariant infinite subgroup of $D^{\gs}_{\infty}=\cup_{n\in\N} \Ker(F^n)\cap\gs$ is dense in
$\gs$.
\end{enumerate}
Then $\pi\mu=\lambda_{\gs}$.
\end{cor}

\section{A discussion about the assumptions}

Comparing the assumptions of Theorems~\ref{G1} and~\ref{G2} with those of Theorems~\ref{HMM1}, \ref{HMM2}
and~\ref{P2} is not completely obvious. Already Theorems~\ref{G1} and~\ref{G2} consider bipermutative algebraic
CA without restriction on the neighborhood. In this section we discuss about the assumptions of these theorems
and show that Theorems~\ref{G1} and~\ref{G2} generalize Theorems~\ref{HMM2} and~\ref{P2} but the ergodic
assumptions cannot be compared with these of Theorem~\ref{HMM1}.

\subsection{Class of CA considered}

Theorems~\ref{G1} and~\ref{G2} consider algebraic bipermutative CA without restriction on the neighborhood. The
bipermutativity is principally used to prove the entropy formula of Lemma~\ref{cG}. We can hope such formula for
expansive CA. Subsection 4.3 gives a result in this direction. The next proposition shows that it is equivalent
to consider algebraic CA or the restriction of a linear CA.

\begin{prop}\label{alglinCA}
Let $\az$ be any Abelian group shift and let $(\az,F)$ be an algebraic CA. There exist $(\bz,G)$ a linear CA and
$\Gamma$ a $\s_{\B^{\Z}}$-invariant subgroup of $\bz$ such that $(\az,\s,F)$ is isomorphic to
$(\Gamma,\s_{\B^{\Z}},G)$ in both dynamical and algebraical sense.
\end{prop}
\begin{proof}
Let $(\az,F)$ be an algebraic CA. By B.P. Kitchens \cite[Proposition 2]{Kitchens-1987}, there exists $\B'$ a
finite Abelian group, $\Gamma$ a Markov subgroup of $\B'^{\Z}$ and $\varphi$ a continuous group isomorphism such
that $\varphi \circ \s = \s_{\B'^{\Z}} \circ  \varphi$. Define $G'=\varphi \circ F \circ \varphi^{-1}$. One has
the following commutative diagram:
\[
\begin{CD}
\az  @>\s , F >>  \az  \\
@VV\varphi V     @VV\varphi V      \\
\Gamma  @>\s_{\B'^{\Z}} , G'>>  \Gamma
\end{CD}
\begin{gathered}
\qquad \qquad \qquad  \varphi \circ  \s =  \s_{\B'^{\Z}} \circ  \varphi \\
\qquad \qquad \qquad  \varphi \circ  F =  G' \circ  \varphi
\end{gathered}
\]

$G'$ is continuous and commutes with $\s_{\B'^{\Z}}$, so it is a CA on $\Gamma'$. We want to extend $G'$ to
obtain a linear CA. By G.A. Hedlund \cite{Hedlund-1969}, there exist a neighborhood $\U$, $H$ a subgroup of
$\B'^{\U}$ and a local function $\overline{G'}:H\to\B'$ which define $G'$. Moreover, by linearity,
$\overline{G'}$ is a group morphism. If we could extend $\overline{G'}$ to a morphism from $\B^{\U}$ to $\B$
(where $\B'$ was a subgroup of $\B$), we would obtain the local rule of a linear CA.

There exist $d,k\in\N$ such that $\B'$ can be viewed as a subgroup of $(\Z/d\Z)^k$. If $\B=(\Z/d\Z)^k$, then $H$
can be viewed as a subgroup of $\B^{\U}$. By the Fundamental Theorem of Finitely Generated Abelian Group
\cite[Theorem 7.8]{Lang-2002}, there exist $e_1,...,e_{k|\U|}$ a basis of $\B^{\U}$ and $a_1,...,a_{k|\U|}\in\N$
such that $\B^{\U}=\bigoplus_{i}\langle e_i\rangle$ and $H=\bigoplus_{i}\langle a_i e_i\rangle$. For all $i\in
[1,k|\U|]$, there exist $f_i\in\B$ such that $\overline{G'}(a_i e_i)=a_i f_i$ because the order of
$\overline{G'}(a_i e_i)$ is at most $\frac{d}{a_i}$. We define the morphism $\overline{G}:\B^{\U}\to\B$ by
$\overline{G}(e_i)=f_i$ for all $i\in[1,k|\U|]$. $\overline{G}$ defines a linear CA on $\bz$ denoted $G$ whose
the restriction is $(\Gamma,G')$.
\end{proof}

\begin{remark}
The study of algebraic CA can be restricted to the study of the restriction of linear CA to Markov subgroups.
\end{remark}

Since we consider $\s$-invariant measures, we can assume that the neighborhood of the CA is $\U=[0,r]$. Moreover
it is easy to show the next Proposition and consider CA of neighborhood $\U=[0,1]$.
\begin{prop}\label{equiv}
Let $(\az,F)$ be a CA of neighborhood $\U=[0,r]$. There is a CA $((\A^r)^{\Z},G)$ of neighborhood $\U=[0,1]$ so
that the topological system $(\az,F)$ is isomorphic to the system $((\A^r)^{\Z},G)$ via the conjugacy
$$\phi_r:(x_i)_{i\in\Z}\in\az\to((x_{[ri,ri+r-1]})_{i\in\Z})\in(\A^r)^{\Z}.$$
Furthermore one has:
\begin{eqnarray*}
(\az,F)\textrm{ is bipermutative} & \Longleftrightarrow &((\A^r)^{\Z},G)\textrm{ is bipermutative};\\
(\az,F)\textrm{ is algebraic} & \Longleftrightarrow & ((\A^r)^{\Z},G)\textrm{ is algebraic};\\
(\az,F)\textrm{ is linear} & \Longleftrightarrow & ((\A^r)^{\Z},G)\textrm{ is linear}.
\end{eqnarray*}
\end{prop}

If $\mu\in\Me(\az)$ is $\s$-totally ergodic then $\phi_r\mu\in\Me((\A^r)^{\Z})$ is $\s_{(\A^r)^{\Z}}$-totally
ergodic. Moreover,  by conjugacy, $h_{\mu}(F)>0$ is equivalent to $h_{\phi_r\mu}(G)>0$. So, as suggested in
\cite{Pivato-2003}, Theorem~\ref{P2} holds for algebraic bipermutative CA without any restriction on the
neighborhood.
\begin{cor}\label{corP2}
Let $\az$ be any Abelian group shift, let $(\az,F)$ be an algebraic bipermutative CA (without restriction on the
neighborhood) and let $\mu$ be an $(F,\s)$-invariant probability measure. Assume that:
\begin{enumerate}
\item $\mu$ is totally ergodic for $\s$;

\item $h_{\mu}(F)>0$;

\item $\Ker(F)$ contains no nontrivial $\s$-invariant subgroups.
\end{enumerate}
Then $\mu=\lambda_{\az}$.
\end{cor}
\begin{remark}
The correspondence holds only if $\mu$ is supposed to be $\s$-totally ergodic. Indeed if $\mu$ is $\s$-ergodic,
$\phi_r\mu$ is not necessarily $\s_{(\A^r)^{\Z}}$-ergodic
\end{remark}

\subsection{Ergodicity of action}

Assumption (1) of Theorem~\ref{G1} characterizes the ergodicity of the action $(F,\s)$ on the measure space
$(\az,\Ba,\mu)$. Since we want to characterize $(F,\s)$-invariant measures, it is natural to assume that $\mu$
is $(F,\s)$-ergodic because every $(F,\s)$-invariant measure can be decomposed into $(F,\s)$-ergodic components.
The next relations are easy to check for an $(F,\s)$-invariant probability measure $\mu$:
$$\mu \textrm{ is } (F,\s)\textrm{-totally ergodic}\  \Rightarrow \ \mu \textrm{ is } \s\textrm{-totally ergodic}\  \Rightarrow \ \mu \textrm{ is } \s\textrm{-ergodic}\ \Rightarrow \ \mu \textrm{ is } (F,\s)\textrm{-ergodic};$$
$$  \mu \textrm{ is } \s\textrm{-totally ergodic}\  \Rightarrow \ \mu \textrm{ is } (F,\s)\textrm{-ergodic and } \I_{\mu}(\s)=\I_{\mu}(\s^{k}) \textrm{ for
every } k\geq 1.$$

Thus, hypothesis (1) of Theorem~\ref{P2} implies hypothesis~(1) and~(2) of Theorem~\ref{G2} which imply
hypothesis~(1) and~(2) of Theorem~\ref{G1}. However, we remark that the ergodicity assumption (1) of
Theorem~\ref{HMM1} cannot be compared with hypothesis of Theorem~\ref{G1}. Indeed, there are probability
measures which are $(F,\s)$-ergodic with $\I_{\mu}(\s)=\I_{\mu}(\s^{k})$ for some $k\geq 1$ which are not
$\s$-ergodic. Conversely there exist probability measures which are $\s$-ergodic with
$\I_{\mu}(\s)\ne\I_{\mu}(\s^{k})$ for some $k\geq 1$.

Secondly, if $\A=\Z/p\Z$ and $F=a\,\id +b\,\s$ on $\az$ then $p-1$ is a multiple of the common period of every
element of $\Ker(F)$. So the spectrum assumption~(2) of Theorem~\ref{HMM2} implies hypothesis~(2) of
Theorems~\ref{G1} and~\ref{G2}. For Theorem~\ref{P2} the total ergodicity of $\mu$ under $\s$ is required. This
property does not seem to be very far from hypothesis~(2) of Theorems~\ref{G1} and~\ref{G2}. But condition (2)
of Theorems~\ref{G1} and~\ref{G2} (concerning the $\s$-invariant set) shows the importance of the algebraic
characteristic of the system. The property of $(F,\s)$-total ergodicity of $\mu$ is more restrictive. With such
an assumption Einsiedler \cite{Einsiedler-2003} proves rigidity results for a class of algebraic actions that
are not necessarily CA. To finish, the next example shows that assumption (2) of Theorems~\ref{G1} and~\ref{G2}
is necessary to obtain the characterization of the uniform Bernoulli measure.

\begin{ex}
Let $\A=\Z/2\Z$ and $F=\id+\s$ on $\az$. We consider the subgroup $X_1=\{x\in\az : x_{2n}=x_{2n+1},\forall
n\in\Z\}$, it is neither $\s$-invariant nor $F$-invariant. Let $X_2=\s(X_1)=\{x\in\az : x_{2n}=x_{2n-1},\forall
n\in\Z\}$, $X_3=F(X_1)=\{x\in\az : x_{2n}=0,\forall n\in\Z\}$ and $X_4=F(X_2)=\{x\in\az : x_{2n+1}=0,\forall
n\in\Z\}$. The set $X=X_1\cup X_2 \cup X_3 \cup X_4$ is $(F,\s)$-invariant. Let $\nu$ be the Haar measure on
$X_1$. We consider $\mu=\frac{1}{4}(\nu+\s\nu+F\nu+F\s\nu)$. It is easy to verify that $\mu$ is an
$(F,\s)$-ergodic measure such that $h_{\mu}(\s)>0$. However $X_i\in\I_{\mu}(\s^2)\smallsetminus\I_{\mu}(\s)$ for
all $i\in[1,4]$, hence hypothesis~(2) is  false, so we cannot apply Theorem~\ref{G1} and $\mu$ it is not the
uniform Bernoulli measure. S. Silberger propose similar constructions in \cite{Silberger-2005}.
\end{ex}

\subsection{Positive entropy}

Corollary \ref{centropie} shows that for a nontrivial bipermutative CA $(\az,F)$, the assumption of positive
entropy of $F$ can be replaced by the positive entropy of $F^n\circ\s^m$ for some $(n,m)\in\N\times\Z$. So the
positive entropy hypothesis~(3) of Theorems~\ref{G1} and~\ref{G2} can be replaced by the positive entropy of the
action $(F,\s)$ in some given direction. We can find this type of assumption in \cite{Einsiedler-2003}.

We can also expect a similar formula for an expansive CA $F$ but in this case we have the inequality:
$h_{\mu}(F)>0$ iff $h_{\mu}(\s)>0$. To begin we show an inequality for a general CA.

\begin{prop}
Let $(\az,F)$ be a CA of neighborhood $\U=[r,s]\ni 0$ (not necessarily the smallest possible one) and $\mu$ be
an $(F,\s)$-invariant probability measure on $\az$. Then $h_{\mu}(F)\leq (r-s)\, h_{\mu}(\s)$.
\end{prop}
\begin{proof}
By definition, for $N\in\N$, $l\in\N$ and $x\in\az$, the knowledge of $x_{[rN-l,sN+l]}$ determines
$(F^n(x)_{[-l,l]})_{n\in[0,N]}$. This means that $\p_{[rN-l,sN+l]}$ is a refinement of $\bigvee_{n=0}^N
F^{-n}(\p_{[-l,l]})$. So for $l\geq\max(s,-r)$ we have:
\begin{eqnarray*}
h_{\mu}(F,\p_{[-l,l]}) & = &  \lim_{N\to\infty} \frac{1}{N} H_{\mu}(\bigvee_{n=0}^N F^{-n}(\p_{[-l,l]})) \\
& \leq & \lim_{N\to\infty}\frac{1}{N} H_{\mu}(\p_{[rN-l,Ns+l]})  \\
& =& \lim_{N\to\infty} -\frac{N(s-r)+2l}{N} \frac{1}{N(s-r)+2l} \sum_{u\in \A^{N(s-r)+2l}} \mu([u])\log(\mu[u])\\
& =& (s-r) \, h_{\mu}(\s).
\end{eqnarray*}
We deduce that $h_{\mu}(F)=\lim_{l\to\infty } h_{\mu}(F,\p_{[-l,l]}) \leq (s-r)\, h_{\mu}(\s)$.
\end{proof}

Let $(\az,F)$ be a positively expansive CA. There exists $r_e$, the {\em constant of expansivity}, such as for
all $x,y\in\az$ if $x\ne y$ there exists $n\in\N$ which verifies $F^n(x)_{[-r_e,r_e]}\ne F^n(y)_{[-r_e,r_e]}$.
Then $(\az,F)$ is topologically conjugate to the one-sided subshift $(S_F,\s)$, where $S_F\subset\bn$, with
$\B=\A^{2r_e+1}$, and where $S_F=\{(F^i(x)_{[-r_e,r_e]})_{i\in\N}:x\in\az\}$, via the conjugacy
$\phi_F:x\in\az\to (F^i(x)_{[-r_e,r_e]})_{i\in\N}\in S_F$. Define $\F:S_F\to S_F$ by
$\F\circ\phi_F(x)=\phi_F\circ\s^{r_e}(x)$ for every $x\in\az$. $(S_F,\F)$ is an invertible one-sided CA. Define
the {\em radius of expansivity} $r_T=\max\{r(\F),r(\F^{-1})\}$.

\begin{prop}
Let $(\az,F)$ be a positively expansive CA and $\mu$ an $(F,\s)$-invariant probability measure, then
$h_{\mu}(F)\geq \frac{1}{r_T} \,h_{\mu}(\s)$.
\end{prop}
\begin{proof}
By definition of $r_T$, for $N\in\N$, $l\geq r_e$ and $x\in\az$, the knowledge of
$(F^n(x)_{[-l,l]})_{n\in[0,r_TN]}$ implies the knowledge of $x_{[-N-l,N+l]}$. This means that
$\bigvee_{n=0}^{r_TN} F^{-n}(\p_{[-l,l]})$ is a refinement of $\p_{[-N-l,N+l]}$. A computation similar to that
in the previous proof shows that $r_T\, h_{\mu}(F)\geq h_{\mu}(\s)$.
\end{proof}

This result can be viewed as a rigidity result. Indeed for an expansive CA $(\az,F)$, the measure entropy of F
and $\s$ are linked for an $(F,\s)$-invariant measure. This is a first step in the research of Lyapunov
exponents for expansive CA \cite{Tisseur-2000}.

\subsection{$(F,\s)$-invariant subgroups of $D_{\infty}$}

Now let us discuss assumption~(4) of Theorems~\ref{G1} and~\ref{G2} which is an algebraic condition on the CA.
We can remark that Theorems~\ref{HMM1} and~\ref{HMM2} have no such assumption because they concern a particular
class of CA which verifies this assumption: $F=a\,\id+b\,\s$ on $(\Z/p\Z)^{\Z}$ with $p$ prime. By
Proposition~\ref{equiv} it is easy to modify the proof of Theorem~\ref{P2} to consider nontrivial algebraic
bipermutative CA without restriction on the neighborhood (Corollary~\ref{corP2}). But it is necessary to compare
the assumption ``$\Ker(F)$ contains no nontrivial $\s$-invariant subgroups'' with ``every $\s$-invariant
infinite subgroup of $D_{\infty}$ is dense in $\az$''. We show that the second property is more general and give
in Subsection 5.1 a general class of examples where it is the case.

If $H\subset\az$, denote by $\langle H \rangle$ the subgroup generated by $H$, $\langle H \rangle_{\s}$ the
smallest $\s$-invariant subgroup which contains $H$ and $\langle H \rangle_{F,\s}$ the smallest
$(F,\s)$-invariant subgroup which contains $H$. Let $\gs$ be a closed $(F,\s)$-invariant subgroup. If
$H\subset\gs$, then we remark that $\langle H \rangle$, $\langle H \rangle_{\s}$ and $\langle H \rangle_{F,\s}$
are subgroups of $\gs$.

\begin{prop}\label{D1}
Let $(\az,F)$ be an algebraic CA and let $\gs$ be a closed $(F,\s)$-invariant subgroup of $\az$. The following
propositions are equivalent:
\begin{enumerate}
\item $\ds_{\infty}$ contains no nontrivial $(F,\s)$-invariant infinite subgroups.

\item There exist $m\in\N$ and $n_0\geq 0$ such that $\ds_{n_0}\subset \langle d\rangle_{F,\s}$ for all $d\in
\partial\ds_{n_0+m}$.

\item There exists $m\in\N$ such that $\ds_{n_0}\subset \langle d\rangle_{F,\s}$ for all $n_0\in\N^{\ast}$ and
$d\in \partial\ds_{n_0+m}$.

\item There exists $m\in\N$ such that $\ds_{1}\subset\langle d\rangle_{F,\s}$ for all $d\in
\partial\ds_{m+1}$.
\end{enumerate}
\end{prop}
\begin{proof}
$(2)\Rightarrow (1)$ Let $\Gamma$ be an $(F,\s)$-invariant infinite subgroup of $\ds_{\infty}$. We prove by
induction that $\ds_n\subset\Gamma$ for all $n\geq n_0$. Since $\Gamma$ is infinite and $\ds_n$ is finite for
all $n\in\N$, we deduce that there exists $n'\geq 0$ such that there exists $d\in\Gamma\cap
\partial\ds_{n'+n_0+m}$. By $F$-invariance of $\Gamma$ we have $F^{n'}(d)\in\Gamma\cap
\partial\ds_{n_0+m}$, thus $\ds_{n_0}\subset \langle F^{n'}(d)
\rangle_{F,\s}\subset\Gamma$.

Let $n\geq n_0$ and assume that $\ds_n\subset\Gamma$. We want to show that $\ds_{n+1}\subset\Gamma$. As before,
since $\Gamma$ is infinite and $F$-invariant we can find $d\in\Gamma\cap \partial\ds_{n+1+m}$. From
$F^{n+1-n_0}(d)\in \Gamma\cap \partial\ds_{n_0+m}$, we deduce $\ds_{n_0}\subset \langle F^{n+1-n_0}(d)
\rangle_{F,\s}$. Let $d'\in \ds_{n+1}$. Then $F^{n+1-n_0}(d')\in \ds_{n_0}\subset \langle F^{n+1-n_0}(d)
\rangle_{F,\s}$ and consequently there exists a finite subset $\V\subset\Z\times\N$ such that
$F^{n+1-n_0}(d')=\sum_{(u,m')\in\V}c_{u,m'}\s^u\circ F^{m'+n+1-n_0}(d)$ where $c_{u,m'}\in\Z$. We deduce that
$d'-\sum_{(u,m')\in\V}c_{u,m'} \s^u\circ F^{m'}(d)\in \ds_{n+1-n_0}\subset \ds_n\subset\Gamma$. But
$d\in\Gamma$, so $\s^n\circ F^{m'}\in\Gamma$ for all $(n,m')\in\V$. Thus, $d'\in\Gamma$. This holds for any
$d'\in\ds_{n+1}$. Thus, $\ds_{n+1}\subset\Gamma$. By induction, $\ds_k\subset\Gamma$ for all $k\in\N$. Finally,
$\ds_{\infty}=\cup_{n\in\N} \ds_n\subset \Gamma$.

$(1)\Rightarrow (4)$ By contradiction, we assume that for all $m \in\N$ there exists $d\in \partial\ds_{m+1}$
such that $\langle d \rangle_{F,\s}\cap \ds_{1} \ne \ds_1$. Since $\ds_1$ is a finite group there exists a
strict subgroup $H$ of $\ds_1$ such that $\Delta=\{d\in \ds_{\infty} | \langle d\rangle_{F,\s}\cap \ds_1 \subset
H\}$ is infinite. Observe that $F(\Delta)\subset\Delta$. For all $d'\in \Delta$ we denote $\Delta_{d'}=\{d\in
\Delta | d'\in\langle d\rangle_{F,\s}\}$. Let $(n_i)_{i\in \N}$ be an increasing sequence such that $\Delta\cap
\partial\ds_{n_i}\ne\emptyset$. If $d\in \Delta \cap \partial\ds_{n_{i+1}}$, we have $d'=F^{n_{i+1}-n_i}(d)\in\langle
d\rangle_{F,\s}$, so that $d\in\Delta_{d'}$, and also $d'\in \Delta\cap \partial\ds_{n_i}$. So we can construct
by induction an infinite sequence $(d_i)_{i\in\N}$ of $\ds_{\infty}$ such that $d_i\in \Delta\cap
\partial\ds_{n_i}$ and $d_{i+1}\in \Delta_{d_i}$ for all $i\in\N$. Thus $\Gamma=\bigcup_{i\in\N} \langle
d_i\rangle_{F,\s}$ is an infinite $(F,\s)$-invariant subgroup of $\ds_{\infty}$ such that $\Gamma\cap
\ds_1\subset H$, which contradicts $(1)$.

$(4)\Rightarrow (3)$ Let $m\in\N$ such that $\ds_{1}\subset\langle d\rangle_{F,\s}$ for all $d\in
\partial\ds_{m+1}$. We prove by induction that for all $n\geq 1$ and $d\in \partial\ds_{n+m}$
one has $\ds_{n}\subset \langle d\rangle_{F,\s}$. For $n=1$ it is the assumption. Assume that the property is
true for $n\in\na$. Let $d\in\partial\ds_{n+1+m}$, since $F^{n}(d)\in\partial\ds_{m+1}$, one has
$\ds_{1}\subset\langle F^{n}(d)\rangle_{F,\s}$. If $d'\in\ds_{n+1}$, then $F^n(d')\in\ds_1$ and we deduce the
existence of $\V\subset\Z\times\N$ such that $F^{n}(d')=\sum_{(u,m')\in\V}c_{u,m'}\s^u\circ F^{m'+n}(d)$ where
$c_{u,m'}\in\Z$. From $d'-\sum_{(u,m')\in\V}c_{u,m'}\s^u\circ F^{m'}(d)\in \ds_n$ and from the fact that
$\ds_n\subset\langle F(d)\rangle_{F,\s}$ because $F(d)\in\partial\ds_{n+m}$, we deduce that
$d'-\sum_{(u,m')\in\V}c_{u,m'}\s^u\circ F^{m'}(d)\in\langle F(d)\rangle_{F,\s}\subset\langle d \rangle_{F,\s}$.
Thus, $d'\in\langle d\rangle_{F,\s}$. One deduces that $\ds_{n+1}\subset\langle d \rangle_{F,\s}$.

$(3)\Rightarrow (2)$ is trivial.
\end{proof}

\begin{cor}\label{ker}
If $\ds_1=\Ker (F)\cap\gs$ contains no nontrivial $\s$-invariant subgroups then $\ds_{\infty}$ contains no
nontrivial $(F,\s)$-invariant infinite subgroups.
\end{cor}
\begin{proof}
If $\ds_1=\Ker (F)\cap\gs$ contains no nontrivial $\s$-invariant subgroups, for all $d\in\partial\ds_1$, the
subgroup $\langle d\rangle_{F,\s}$ must be equal to $\ds_1$. By Proposition~\ref{D1}, one deduce that
$\ds_{\infty}$ contains no nontrivial $(F,\s)$-invariant infinite subgroups.
\end{proof}

For a linear CA $(\az,F)$ where $\A=\Z/n\Z$, the $\s$-invariant subgroups coincide with the $(F,\s)$-invariant
subgroups. From Corollary~\ref{ker} we get directly that Theorem~\ref{G1} is stronger than Theorem~\ref{P2} in
this case. Moreover, if we consider the case of the Theorem~\ref{HMM2}, that is to say that $\A=\Z/p\Z$ with $p$
prime and $F=a\,\id + b\,\s$ with $a\ne 0$ and $b\ne 0$, then $\Ker(F)\simeq\Z/p\Z$ does not contain nontrivial
$\s$-invariant subgroups. So Theorem~\ref{G1} generalizes also Theorem~\ref{HMM2}.

When $\A$ is not cyclic, the $\s$-invariant subgroups does not necessarily coincide with the $(F,\s)$-invariant
subgroups. In this case we do not know if Theorem~\ref{G1} implies Theorem~\ref{P2}. However Corollary~\ref{ker}
implies that Theorem~\ref{G2} is stronger than Theorem~\ref{P2} for every algebraic bipermutative CA.

\section{Extensions to some linear CA}

\subsection{The case $\A=\Z/p\Z$}

Theorem \ref{HMM2} concerns linear CA on $(\Z/p\Z)^{\Z}$ of smallest neighborhood $\U=[0,1]$. We will show that
this implies the fourth assumption of Theorem \ref{G1}. In fact we can show that the fourth assumption is
directly implied when we consider a nontrivial linear CA on $(\Z/p\Z)^{\Z}$. This allows us to prove the
following result.

\begin{prop}\label{czpz}
Let $\A=\Z/p\Z$, let $(\az,F)$ be a nontrivial linear CA with $p$ prime and let $\mu$ be an $(F,\s)$-invariant
probability measure on $\az$. Assume that:
\begin{enumerate}
\item $\mu$ is ergodic for the $\N\times\Z$-action $(F,\s)$;

\item $\I_{\mu}(\s)=\I_{\mu}(\s^{p\,p_1})$ with $k\in\N^{\ast}$ and $p_1$ the smallest common period of all
elements of $\Ker(F)$;

\item $h_{\mu}(F)>0$.
\end{enumerate}
Then: (a) $\mu=\lambda_{\az}$.

(b) Moreover $p_1$ divides $\prod_{i=0}^{r-1}(p^r-p^i)$ where $r=\max\{\U,0\}-\min\{\U,0\}$ and $\U$ is the
smallest neighborhood of~$F$.
\end{prop}
\begin{proof}
Proof of (a): By $(F,\s)$-invariance of $\mu$, we can compose $F$ with $\s$ and assume that the smallest
neighborhood of $F$ is $[0,r]$ with $r\in\N\setminus\{0\}$. So $F=\sum_{u\in [0,r]}f_u\circ\s^u=P_F(\s)$ where
$P_F$ is a polynomial with coefficients in $\Z/p\Z$ with $f_0\ne 0$ and $f_r\ne 0$. We remark that $F$ is
bipermutative.

Case 1: First we assume that $P_F$ is irreducible on $\Z/p\Z$. We can view $D_1(F)$ as a $\Z/p\Z$ vector space
and consider the isomorphism $\s_1:D_1(F)\to D_1(F)$, the restriction of $\s$ at the subgroup $D_1(F)$. By
bipermutativity of $F$, $D_1\simeq(\Z/p\Z)^r$. Moreover $P_F(\s_1)=0$; since $P_F$ is irreducible and its degree
is equal to the dimension of $D_1$, we deduce that $P_F$ is the characteristic polynomial of $\s_1$. Since $P_F$
is irreducible, $D_1(F)$ is $\s_1$-simple, so $D_1(F)$ contains no nontrivial $\s$-invariant subgroups,
see~\cite[\S VI.8]{Arnaudies-Bertin-1993} for more detail. By Corollary~\ref{ker}, $D_{\infty}(F)$ also contains
no nontrivial $(F,\s)$-invariant infinite subgroup, so hypothesis~(4) of Theorem~\ref{G1} is verified.

Case 2: Now we assume that $P_F=P^{\alpha}$ where $P$ is irreducible on $\Z/p\Z$ and $\alpha\in\N$. We have
$D_n(P_F(\s))=\Ker(P^{\alpha n}(\s))=D_{\alpha n}(P(\s))$ for all $n\in\N$. So
$D_{\infty}(P_F(\s))=D_{\infty}(P(\s))$. Now we are in the previous case and the fourth condition of Theorem
\ref{G1} is verified.

Case 3: In the general case $P_F=P_1^{\alpha_1}...P_l^{\alpha_l}$ where $P_i$ is irreducible and $\alpha_i\in\N$
for all $i\in [1,l]$. Let $\Gamma$ be an $(F,\s)$-invariant infinite subgroup of $D_{\infty}(P_F(\s))$. By the
kernel decomposition Lemma~\cite[\S VI.4]{Arnaudies-Bertin-1993}, we have
$D_n(P_F(\s))=D_n(P_1^{\alpha_1}(\s))\oplus...\oplus D_n(P_l^{\alpha_l}(\s))$ for every $n\in\N$. Moreover
$D_n(P_F(\s))\cap\Gamma$ is a $\s$-invariant subspace of $D_n(P_F(\s))$ considered as a $\Z/p\Z$-vector space
and $D_n(P_F(\s))\cap\Gamma=(D_n(P_1^{\alpha_1}(\s))\cap\Gamma)\oplus...\oplus
(D_n(P_l^{\alpha_l}(\s))\cap\Gamma)$. We deduce that
$$D_{\infty}(P_F(\s))\cap\Gamma=\bigoplus_{i\in[1,l]}(D_{\infty}(P_i^{\alpha_1}(\s))\cap\Gamma)\eg{(*)} \bigoplus_{i\in[1,l]} (D_{\infty}(P_i(\s))\cap\Gamma),$$
where $(*)$ follows as in Case 2. There exists $i\in[1,l]$ such that $\Gamma\cap D_{\infty}(P_i(\s))$ is an
infinite subgroup. By Case 1, one has $\Gamma\cap D_{\infty}(P_i(\s))=D_{\infty}(P_i(\s))$, so
$D_{\infty}(P_i(\s))\subset \Gamma$. We deduce that $\Gamma$ is dense, because $D_{\infty}(P_i(\s))$ is dense,
because $P_i(\s)$ is bipermutative. Thus the fourth condition of Theorem~\ref{G1} is verified; part (a) of the
proposition follows.

Proof of (b): If $x\in\Ker(F)$, then the coordinates of $x$ verify
$x_{n+r}=-f_r^{-1}\sum_{i=0}^{r-1}f_{i}x_{n+i}$ for all $n\in\Z$. This recurrence relation can be expressed with
a matrix. For all $n\in\Z$ one has $X_{n+1}=AX_n$ where
$$X_n=\begin{pmatrix}x_{n+r-1}\\\vdots\\x_n\end{pmatrix}\textrm{ and }
A=\left[\begin{matrix} -f_{r-1}f_r^{-1}&\cdots&\cdots&\cdots& -f_{0}f_r^{-1}\\
1&0&\cdots&\cdots&0\\
0&\ddots&\ddots&&\vdots\\
\vdots&\ddots&\ddots&\ddots&\vdots\\
0&\cdots&0&1&0\\
\end{matrix}\right].$$
$A$ is invertible because $f_0\ne 0\ne f_r$, and for all $n\in\Z$ one has $X_n=A^nX_0$. Thus the period of $X_n$
divides the period of $A$, which divides the cardinality of the set of invertible matrices on $\Z/p\Z$ of size
$r$, that is to say the number of bases of $(\Z/p\Z)^r$, which is $\prod_{i=0}^{r-1}(p^r-p^i)$.
\end{proof}
\begin{remark}
Proposition~\ref{czpz} still holds if $((\Z/p\Z)^{\Z},F)$ is an affine CA.
\end{remark}

\begin{remark}
Proposition~\ref{czpz} extends to the case when $\A$ is a finite field and $F=\sum_{u\in\U}f_u\s^u$ is a linear
CA where each coefficient $f_u$ is the multiplication by an element of the field.
\end{remark}

Let $((\Z/p\Z)^{\Z},F)$ be a nontrivial linear CA where $P_F(\s)=\sum_{u\in [0,r]}f_u\circ\s^u$ is a polynomial
with coefficients in $\Z/p\Z$ with $f_0\ne 0$ and $f_r\ne 0$. In this case Theorem \ref{P2}, generalized to
nontrivial algebraic bipermutative CA without restriction on the neighborhood, holds only if $\Ker(F)$ contains
no nontrivial $\s$-invariant subgroups, which is equivalent to the irreducibility of $P_F$.
Proposition~\ref{czpz} holds for every linear CA on $(\Z/p\Z)^{\Z}$.

\subsection{The case $\A=\Z/p\Z\times\Z/q\Z$}

Now we consider $\A=\Z/p\Z\times \Z/q\Z$ with $p$ and $q$ distinct primes and $(\az,F)$ a linear bipermutative
CA. In this case $D_{\infty}$ contains infinite $\s$-invariant subgroups which are not dense in $\az$. For
example $D^{\Gamma_1}_{\infty}$ and $D^{\Gamma_2}_{\infty}$ where $\Gamma_1=(\Z/p\Z)^{\Z}\times
\{0_{(\Z/q\Z)^{\Z}}\}$ and $\Gamma_2=\{0_{(\Z/p\Z)^{\Z}}\}\times (\Z/q\Z)^{\Z}$. The measures
$\lambda_{\Gamma_1}$ and $\lambda_{\Gamma_2}$ are $(F,\s)$-totally ergodic with positive entropy for $\s$. If
$\mu$ is an $(F,\s)$-invariant measure which verifies conditions of the Theorem~\ref{G1}, we cannot conclude
that $\mu=\lambda_{\az}$. But if we consider the natural factor $\pi_1:\az\to\Gamma_1$ and
$\pi_2:\az\to\Gamma_2$, then by Corollary~\ref{cfactor}, one has $\pi_1\mu=\lambda_{\Gamma_1}$ or
$\pi_2\mu=\lambda_{\Gamma_2}$. A natural conjecture is this: if every cellular automaton factor of $F$ has
positive entropy, then $\mu=\lambda_{\az}$. The problem is to rebuild the measure starting from $\pi_1\mu$ and
$\pi_2\mu$.

\subsection{The case $\A=\Z/p^k\Z$}

In this case we do not know under what extra conditions an $(F,\s)$-invariant measure is the Haar measure.
Moreover some linear CA are not bipermutative. The next lemma shows how to remove this condition when you
consider a power of the CA.

\begin{lemma}\label{zpkzlemma}
Let $(\az,F)$ be a linear CA with $\A=\Z/p^k\Z$, where $p$ is prime, $k\geq 1$ and $F=\sum_{i\in [r,s]}f_i\,
\s^i$, with $f_i\in\Z/p^k\Z$. Let $\widehat{\U}=\{i\in[r,s]:f_i\textrm{ coprime with }p\}$,
$\hat{r}=\min\widehat{\U}$ and $\hat{s}=\max\widehat{\U}$. Assume $\widehat{\U}$ is not empty and
$\hat{r}<\hat{s}$.

Then $F^{p^{k-1}}$ is bipermutative of smallest neighborhood $\U'=[p^{k-1}\hat{r},p^{k-1}\hat{s}]$.
\end{lemma}
\begin{proof}
We can write $F=P_F(\s)$ with $P_F\in\Z/p^k\Z[X,X^{-1}]$. We decompose $P_F=P_1+pP_2$ where
$P_1=\sum_{i\in\widehat{\U}}f_i\,X^i$. By Fermat's little theorem and induction on $j\geq 1$, we can easily
prove that:
$$(P_1+pP_2)^{p^j}= (P_1)^{p^j} \mod p^{j+1}.$$
So we have $P_F^{p^{k-1}}=P_1^{p^{k-1}}=\sum_{i\in[p^{k-1}\hat{r},p^{k-1}\hat{s}]} g_i X^i$ where
$g_i\in\Z/p^k\Z$. Moreover $g_{p^{k-1}\hat{r}}=f_{\hat{r}}^{p^{k-1}}$ and
$g_{p^{k-1}\hat{s}}=f_{\hat{s}}^{p^{k-1}}$ are relatively prime to $p$. We deduce that
$F^{p^{k-1}}=P_F^{p^{k-1}}(\s)$ is bipermutative of smallest neighborhood $\U'=[p^{k-1}\hat{r},p^{k-1}\hat{s}]$.
\end{proof}

Now we can deduce from Corollary \ref{centropie} an entropy formula for general linear CA on $(\Z/p^k\Z)^{\Z}$.
\begin{cor}
Let $(\az,F)$ be a linear CA with $\A=\Z/p^k\Z$, where $p$ is prime, $k\geq 1$, and $F=\sum_{i\in [s,r]}f_i\,
\s^i$ with $f_i\in\Z/p^k\Z$. Let $\hat{r}<\hat{s}$ be as in Lemma~\ref{zpkzlemma}. Let $\mu$ be an
$(F,\s)$-invariant probability measure on $\az$. Then $h_{\mu}(F)=(\max(\hat{r},0)-\min(\hat{s},0))h_{\mu}(\s)$.
\end{cor}

\begin{cor}\label{corzpk}
Let $(\az,F)$ be a linear CA with $\A=\Z/p^k\Z$, where $p$ is prime, $k\geq 1$, and $F=\sum_{i\in [s,r]}f_i\,
\s^i$ with $f_i\in\Z/p^k\Z$. Assume that for at least two $i\in [s,r]$, $f_i$ is relatively prime with $p$. Let
$\gs$ be a closed $(F,\s)$-invariant subgroup of $\az$ and let $\mu$ be an $(F,\s)$-invariant probability
measure on $\az$ with $\supp(\mu)\subset\gs$. Assume that:
\begin{enumerate}
\item $\mu$ is ergodic for the $\N\times\Z$-action induced by $(F,\s)$;

\item $\I_{\mu}(\s)=\I_{\mu}(\s^{p p_1})$ with $p_1$ the smallest common period of all elements of $\Ker(F)$;

\item $h_{\mu}(\s)>0$;

\item every $\s$-invariant infinite subgroup of $\ds_{\infty}(F)=\cup_{n\in\N} \Ker(F^n)\cap\gs$ is dense in
$\gs$.
\end{enumerate}
Then $\mu=\lambda_{\gs}$.
\end{cor}

\begin{ex}
Let $\A=\Z/4\Z$, we consider the CA $(\az,F)$ defined by $F=\id+\s+2\s^2$. Then $\gs=\{0,2\}^{\Z}$ satisfies the
conditions of Corollary \ref{corzpk}. In this case the only $(F,\s)$-invariant probability measure of positive
entropy known are $\lambda_{\az}$ and $\lambda_{\gs}$.
\end{ex}

\section{Measure rigidity for some affine one-sided expansive CA}

An invertible onesided CA $(\an,F)$ is called expansive if there exists a constant $r_e\in\N$ such that for all
$x,y\in\an$, if $x\ne y$ there exists $n\in\Z$ which verifies $F^n(x)_{[0,r_e]}\ne F^n(y)_{[0,r_e]}$. Expansive
CA are different from positively expansive CA because we look also the past of the orbit. M. Boyle and A. Maass
introduced in \cite{Boyle-Maass-2000} a class of onesided invertible expansive CA which have remarkable
combinatorial properties. Further properties were obtained in \cite{Dartnell-Maass-Schwartz-2003}. We study this
class of examples from the point of view of measure rigidity. This class of CA is not bipermutative so we cannot
apply directly Theorem~\ref{G1}. However, in some case, it is possible to associate a ``dual" CA which
correspond to the assumptions of Theorem~\ref{G1}. This is a first step to study measure rigidity for expansive
CA.

We are going to recall some properties obtained in \cite{Boyle-Maass-2000}. Let $F:\an\to\an$ be a CA such that
$r(F)=1$. Associate to $F$ the equivalence relation over $\A$: $a\re_F b$ iff
$\overline{F}(\cdot\,a)=\overline{F}(\cdot\,b)$ as a function from $\A$ to $\A$; and we write $\p_{\re_F}$ the
partition induced by $\re_F$ and $C_{\re_F}(a)$ the class associated to $a$. Define also $\pi_F:\A\to\A$ by
$\pi_F(a)=\overline{F}(aa)$ for any $a\in\A$.

\begin{prop}[\cite{Boyle-Maass-2000}]\label{BM1}
A onesided  CA $F:\an\to\an$ with $r(F)=1$ is invertible with $r(F^{-1})=1$ iff the following conditions hold:
\begin{enumerate}
\item $\pi_F$ is a permutation,

\item $F$ is left permutative,

\item $\forall a\in A$, $\Succ_F(a):=\im(\overline{F}(a\,\cdot))\subset \pi_F(C_{\re_F}(a))$.
\end{enumerate}
\end{prop}

If $F$ is an expansive invertible CA with $r(F)=r(F^{-1})=1$, then $(\an,F)$ is topologically conjugate to the
bilateral subshift $(S_F,\s)$ where  $S_F=\{(F^i(x)_0)_{i\in\Z}:x\in\an\}$ via the conjugacy
$\phi_F:x\in\an\to(F^i(x)_0)_{i\in\N}\in S_F$. Define $\F:S_F\to S_F$ by $\F(\phi_F(x))=\phi_F(\s(x))$ for every
$x\in\an$. If $F$ is expansive then $(S_F,\F)$ is a CA (defined on $S_F$ instead of a fullshift). Invertible
expansive CA with $r(F)=r(F^{-1})=r(\F)=1$ can be characterized as follows:

\begin{prop}[\cite{Boyle-Maass-2000}]\label{BM2}
A onesided invertible CA $F:\an\to\an$ with $r(F)=r(F^{-1})=1$ is expansive with $r(\F)=1$ iff the following
conditions are verified:
\begin{enumerate}
\item $|C\cap \pi_F(C')|\leq 1$ for any $C$, $C'\in\p_{\re_F}$,

\item $\forall a\in\A$, $\Succ_F(a):=\im(\overline{F}(a\,\cdot))=\pi_F(C_{\re_F}(a))$.
\end{enumerate}
\end{prop}
Such a CA is said to be in Class (A). The alphabet $\A$ of a CA in Class (A) has cardinality $n^2$ for some
$n\in\N$.

Write $\B=\p_{\re_F}$. In \cite{Boyle-Maass-2000}, the authors show that $(S_F,\s)$ is conjugate to the full
shift $(\bz,\s)$ by $\varphi :S_F\to \bz$ such that $\varphi((a_i)_{i\in\Z})=(C_{\re_F}(a_i))_{i\in\Z}$. The CA
$(S_F,\F)$ determines by $\varphi$ a CA $(\bz,\Ft)$ on $\bz$ and $(S_F,\F)$ is conjugate to $(\bz,\Ft)$. To sum
up we have:
\[
\begin{array}{ccccc}
(\an,\s) & \equiv & (S_F,\F) & \equiv &(\bz,\Ft),\\
(\an,F) & \equiv & (S_F,\s) & \equiv & (\bz,\s),
\end{array}
\]
(where $\equiv$ means topologically conjugate).

\begin{prop}
If $F$ is in Class (A) then $\Ft$ is bipermutative.
\end{prop}
\begin{proof}
Let $(\an,F)$ be a CA in the class (A) and let $\alpha,\alpha',\beta,\gamma,\delta\in\B$ such that
$\overline{\Ft}(\alpha,\beta,\gamma)=\overline{\Ft}(\alpha',\beta,\gamma)=\delta$. Suppose $\beta=\varphi(b)$,
for some $b\in S_F$. Then $b\in\pi_F(\gamma)$ by condition~(2) of Proposition~\ref{BM2}, so
$b\in\beta\cap\pi_F(\gamma)$, which is a singleton set by condition~(1). Hence $\beta$ and $\gamma$ uniquely
determine $b$. Likewise, if $\alpha=\varphi(a)$ and $\alpha'=\varphi(a')$ for some $a,a'\in S_F$, then we must
have $a,a'\in\overline{F}(b,\delta)$. But $\overline{F}(b,.):\A\to\A$ is constant on $\delta$ by definition of
the partition $\p_{\re_F}$, so $a=a'$ so $\alpha=\alpha'$. We deduce that the function
$\overline{\Ft}(\cdot,\beta,\gamma):\B\to\B$ is injective. So it is bijective because $\B$ is finite. Thus,
$(\B,\Ft)$ is left-permutative.

In the same way we can prove that $(\B,\Ft)$ is right-permutative by applying Propositions~\ref{BM1}
and~\ref{BM2} to $F^{-1}$ instead. The result follows.
\end{proof}
A natural question after this proposition is to characterize the CA $F$ in class (A) such that $\Ft$ is
algebraic to apply previous theorems. We have only the next sufficient condition:
\begin{prop}\label{classA}
Let $(\an,F)$ be a linear CA, $F=f_0\id+f_1\s$ where $f_0$ and $f_1$ are endomorphisms of $\A$ extended
coordinate by coordinate to $\an$.

(a) $F$ is invertible with $r(F^{-1})=1$ iff $f_0$ is an automorphism and $f_1\circ f_0^{-1}\circ f_1=0$.

(b) $F$ is in Class (A) iff $f_0$ is an automorphism, $\im f_1=f_0(\Ker f_1 )$ and $\im f_1\cap\Ker f_1=\{0\}$.

(c) When $(\az,F)$ is in Class (A), the CA $(\p_{\re_F}^{\Z},\Ft)$ is linear.
\end{prop}
\begin{proof}
First we remark that $b\in C_{\re_F}(b')$ iff $f_0(a)+f_1(b)=f_0(a)+f_1(b')$ for all $a\in\A$; this is
equivalent to $b\in b'+\Ker f_1 $. So $C_{\re_F}(b)=b+\Ker f_1 $ for all $b\in\A$. Thus, $\p_{\re_F}\cong
\A/\Ker f_1 $. Moreover $\Succ_F(a)=\im(\overline{F}(a\,\cdot))=f_0(a)+\im f_1$ for all $a\in\A$, and
$\pi_F=f_0+f_1$.

Proof of (a): Assuming $f_0$ is an automorphism and $f_1\circ f_0^{-1}\circ f_1=0$, it is possible to express
$F^{-1}$ as: $F^{-1}=f_0^{-1}\id - f_0^{-1}\circ f_1\circ f_0^{-1}\s$. Conversely, if $F$ is invertible with
$r(F^{-1})=1$, by Proposition~\ref{BM1}, $f_0$ is an automorphism because $F$ is left-permutative and $f_1\circ
f_0^{-1}\circ f_1=0$ because for some $a\in\A$ one has:
$$f_0(a)+\im (f_1)=\Succ_F(a)\subset\pi_F(C_{\re_F}(a))=f_0(a)+f_1(a)+f_0(\Ker f_1 ),$$
 that is to say $\im f_1\subset f_0(\Ker f_1 )$.

Proof of (b): As in the proof of (a), one has $\Succ_F(a)=\pi_F(C_{\re_F}(a))$ for any $a\in\A$ iff $\im
f_1=f_0(\Ker f_1)$. Moreover, if $|C\cap\pi_F(C')|\leq 1$ for any $C,C'\in\p_{\re_F}$, then $0+\Ker f_1\cap
\pi_F(0+\Ker f_1)=\Ker f_1\cap f_0(\Ker f_1)=\Ker f_1\cap\im f_1 =\{0\}$. Conversely, for any $b,b'\in \A$ one
has $C_{\re_F}(b)\cap\pi_F(C_{\re_F}(b))=b+\Ker f_1\cap \pi_F(b')+\im f_1$, so if $\Ker f_1\cap\im f_1 =\{0\}$
then $C_{\re_F}(b)\cap\pi_F(C_{\re_F}(b))$ contains at most one element. Characterization of linear CA in Class
(A) follows from Proposition~\ref{BM2}.

Proof of (c): Let $(\an,F)$ be a CA in the class (A). We will show that $(\p_{\re_F}^{\Z},\Ft)$ is linear. Since
$\A$ is finite Abelian and $\im f_1\cap\Ker f_1=\{0\}$ by (b), one has $\im f_1\oplus\Ker f_1=\A$. Moreover $\im
f_1$ and $\Ker f_1$ are isomorphic to the same group, denoted $\B$, because $f_0$ is an automorphism and $\im
f_1=f_0(\Ker f_1 )$ by (b). An element $a\in\A$ is written $\binom{x}{y}$ where $x\in\im f_1\simeq \B$ and
$y\in\Ker f_1\simeq\B$. One has $\p_{\re_F}\simeq\A/\Ker f_1 \simeq\im f_1\simeq\B$. We want to show that
$(\B^{\Z},\Ft)$ is linear. We can write $f_0$ and $f_1$ as $2\times 2$-matrices with coefficients in $\Hom(\B)$:
\[
f_0= \left[\begin{matrix} f_{0,11}&f_{0,12}\\f_{0,21}&f_{0,22}
 \end{matrix}\right]
 \textrm{ and } f_1=\left[\begin{matrix} f_{1,11}&f_{1,12}\\f_{1,21}&f_{1,22}
 \end{matrix}\right].
\]

Since $\im f_1=f_0(\Ker f_1)$ one has $f_{0,22}=0$ and since $f_0$ is an automorphism we deduce that $f_{0,12}$
and $f_{0,21}$ are automorphisms of $\B$. Since the second coordinate corresponds to the kernel of $f_1$, one
has $f^1_{12}=f^1_{22}=0$ and since $\im f_1\cap\Ker f_1=\{0\}$ one has $f_{1,21}=0$. Moreover $f_{1,11}$ is an
automorphism of $\B$ since it is the restriction of $f_1$ at $\im f_1$. So we have:
\[
f_0= \left[\begin{matrix} f_{0,11}&f_{0,12}\\f_{0,21}& 0
 \end{matrix}\right],
f_0^{-1}= \left[\begin{matrix} 0&f^{-1}_{0,21}\\f^{-1}_{0,12}& - f^{-1}_{0,12}\circ f_{0,11}\circ f^{-1}_{0,21}
 \end{matrix}\right]
 \textrm{ and } f_1=\left[\begin{matrix} f_{1,11}&0\\0&0
 \end{matrix}\right].
\]
These formulas are illustrated by the next diagram which represents the action of $\overline{F}$ and
$\overline{F^{-1}}$ on a neighborhood:
\[
\begin{array}{cccccc}
\vdots\\
\begin{pmatrix} f_{0,11}(x_0)+f_{0,12}(y_0)+f_{1,11}(x_1) \\ f_{0,21}(x_0) \end{pmatrix}\\
\\
\begin{pmatrix} x_0\\ y_0 \end{pmatrix} & \begin{pmatrix} x_1 \\ y_1  \end{pmatrix}  & \cdots     \\
\\
\begin{pmatrix} f^{-1}_{0,21}(y_0) \\ f^{-1}_{0,12}(x_0) - f^{-1}_{0,12}\circ f_{0,11}\circ f^{-1}_{0,21}(y_0)+ f^{-1}_{0,12}\circ f_{1,11}\circ f^{-1}_{0,21}(y_1)\end{pmatrix}\\
\vdots
\end{array}
\]
We deduce that $\Ft= f^{-1}_{1,11}\circ\s - f^{-1}_{1,11}\circ f_{0,11}\circ\id - f^{-1}_{1,11}\circ
f_{0,12}\circ f_{0,21} \circ\s^{-1}$, so $(\bz,\Ft)$ is linear.
\end{proof}
With Proposition \ref{czpz} and the conjugacy relations it is possible to characterize the uniform Bernoulli
measure of some linear CA in Class (A):
\begin{prop}\label{rigClassA}
Let $(\an,F)$ be an affine invertible CA in Class (A) with $|\A|=p^2$ with $p$ prime. Let $\mu$ be an
$(F,\s)$-invariant probability measure on $\an$. Assume that:
\begin{enumerate}
\item $\mu$ is ergodic for the $\Z\times\N$-action $(F,\s)$;

\item $\I_{\mu}(F)=\I_{\mu}(F^{p(p-1)(p^2-1)})$;

\item $h_{\mu}(\s)>0$.
\end{enumerate}
Then $\mu=\lambda_{\an}$.
\end{prop}
\begin{proof}
By Proposition~\ref{classA}, $\Ft$ is a linear bipermutative CA of neighborhood $[-1,1]$ on $\bz$, where
$\B=\Z/p\Z$ . There exist $\phi:\an\to\bz$ such that $(\an,F,\s)$ and $(\bz,\s,\Ft)$ are conjugate via $\phi$,
so:
\begin{enumerate}
\item $\phi\mu$ is ergodic for the $\N\times\Z$-action $(\Ft,\s)$;

\item $\I_{\phi\mu}(\s)=\I_{\phi\mu}(\s^{p(p-1)(p^2-1)})\eg{(\ast)}\I_{\phi\mu}(\s^{p(p-1)})$, where $(\ast)$ is
by Remark~\ref{annexelemma};

\item $h_{\phi\mu}(\Ft)>0$.
\end{enumerate}
By Proposition~\ref{czpz}(a) we deduce that $\phi\mu=\lambda_{\bz}$ so $\mu=\lambda_{\an}$.
\end{proof}

The next example shows two CA of class (A) with $|\A|=2^2$.

\begin{ex}
Let $\A=\Z/2Z\times\Z/2\Z$, we define two CA $(\an,F_1)$ and $(\an,F_2)$ in Class (A) by:
\begin{gather*}
\overline{F_1}\left(\begin{pmatrix}x_0 \\ y_0 \end{pmatrix}\begin{pmatrix}x_1 \\ y_1 \end{pmatrix}\right) =
\left[\begin{matrix} 1 & 1 \\ 1 & 0 \end{matrix}\right] \begin{pmatrix}x_0 \\ y_0 \end{pmatrix} +
\left[\begin{matrix} 1 & 0 \\ 0 & 0 \end{matrix}\right] \begin{pmatrix}x_1 \\ y_1 \end{pmatrix}\\
\textrm{ and } \\
\overline{F_2}\left(\begin{pmatrix}x_0 \\ y_0 \end{pmatrix}\begin{pmatrix}x_1 \\ y_1 \end{pmatrix}\right) =
\left[\begin{matrix} 0 & 1 \\ 1 & 0 \end{matrix}\right] \begin{pmatrix}x_0 \\ y_0 \end{pmatrix} +
\left[\begin{matrix} 1 & 0 \\ 0 & 0 \end{matrix}\right] \begin{pmatrix}x_1 \\ y_1 \end{pmatrix}\\
\end{gather*}

The first coordinate corresponds to the class of $\p_{\re_{F_i}}$ and the second coordinate corresponds to the
class of $\p_{\re_{F_i^{-1}}}$. For $i\in\{1,2\}$, let $\mu_i$ be such that:
\begin{enumerate}
\item $\mu_i$ is $(F_i,\s)$-ergodic and $\I_{\mu_i}(F)=\I_{\mu_i}(F^{6})$

\item $\exists (n,m)\in\N\times\Z$ such that $h_{\mu_i}(\s^n\circ F_i^m)>0$
\end{enumerate}

All the hypothesis of Proposition~\ref{rigClassA} are satisfied, we can conclude that $\mu_i=\lambda_{\an}$ for
all $i\in\{1,2\}$. To see where Theorem~\ref{P2} does not hold when we assume $\mu$ $\s$-totally ergodic, we are
going to exhibit $\Ker (\widetilde{F^T_i})$ for $i\in\{1,2\}$.

For $F_1$ one has:
$$\overline{\widetilde{F_1^T}}(\alpha,\beta,\gamma)=\alpha + \beta + \gamma.$$

This formula is illustrated by the next diagram which represents the action of $\overline{F_1}$ and
$\overline{F_1^{-1}}$ on a neighborhood:
\[
\begin{array}{cccccc}
{ \bf \vdots}\\
{\bf \begin{pmatrix} x_0+x_1+y_0 \\ x_0 \end{pmatrix}}\\
\\
\mathbf{\begin{pmatrix} x_0\\ y_0 \end{pmatrix}} & \mathbf{\begin{pmatrix} x_1 \\ y_1  \end{pmatrix}}  & \mathbf{\cdots}     \\
\\
\mathbf{\begin{pmatrix} y_0 \\ y_0+y_1+x_0 \end{pmatrix}}\\
\mathbf{\vdots}
\end{array}
\]

So we have:
$$D_1(\widetilde{F^T_1})= \Ker (\widetilde{F^T_1}) =  \{ ^{\infty}000^{\infty} , ^{\infty}011^{\infty} , ^{\infty}110^{\infty} , ^{\infty}101^{\infty} \} \cong \Z/2\Z\times\Z/2\Z.$$
$\Ker(F^T_1)$ contains no nontrivial $\s$-invariant subgroups. Then $\mu_1= \lambda_{\an}$ by Theorem~\ref{G1}
and Corollary~\ref{ker}. In this case, if $\mu$ was $\s$-totally ergodic, then we could have also applied
Theorem~\ref{P2} to conclude that $\mu=\lambda_{\an}$.

For $F_2$ one has:
$$\overline{\widetilde{F_2^T}}(\alpha,\beta,\gamma)=\alpha + \gamma.$$

This formula is illustrated by the next diagram which represents the action of $\overline{F_2}$ and
$\overline{F_2^{-1}}$ on a neighborhood:
\[
\begin{array}{cccccc}
\vdots\\
\begin{pmatrix} x_1+y_0 \\ x_0 \end{pmatrix}\\
\\
\begin{pmatrix} x_0\\ y_0 \end{pmatrix} & \begin{pmatrix} x_1 \\ y_1  \end{pmatrix}  & \cdots     \\
\\
\begin{pmatrix} y_0 \\ y_1+x_0 \end{pmatrix}\\
\vdots
\end{array}
\]

One obtains:
$$D_1(\widetilde{F^T_2})= \Ker (\widetilde{F^T_2}) =  \{ ^{\infty}00^{\infty} , ^{\infty}11^{\infty} , ^{\infty}01^{\infty} , ^{\infty}10^{\infty} \} \cong
\Z/2\Z\times\Z/2\Z, $$
$$ D_2(\widetilde{F^T_2}) =  \langle D_1 \cup \{^{\infty}0001^{\infty} , ^{\infty}0111^{\infty} , ^{\infty}0011^{\infty}
\}\rangle_{\s}.
$$
We remark that $\forall d\in \partial D_{2}$ one has $D_1\subset\langle d\rangle_{\s}$, so $\mu_2=\lambda_{\an}$
by Theorem~\ref{G1} and Proposition~\ref{D1}. We can also remark that in this case $\{ ^{\infty}00^{\infty} ,
^{\infty}11^{\infty}\}$ is a nontrivial $\s$-invariant subgroup of $\Ker(F)$ so Theorem~\ref{P2} would not
apply, even if we assumed that $\mu$ was $\s$-totally ergodic.
\end{ex}

\paragraph*{Acknowledgments:}
I would like to thank Alejandro Maass and Francois Blanchard for many stimulating conversations during the
process of writing. I thank also the referee, the final version of this paper owe to his numerous suggestions. I
thank Nucleus Millennium P01-005 and Ecos-Conicyt C03E03 for financial support. Last but not least special
acknowledgment to Wen Huang for his attentive reading and his relevant remarks.

\def\ocirc#1{\ifmmode\setbox0=\hbox{$#1$}\dimen0=\ht0 \advance\dimen0
  by1pt\rlap{\hbox to\wd0{\hss\raise\dimen0
  \hbox{\hskip.2em$\scriptscriptstyle\circ$}\hss}}#1\else {\accent"17 #1}\fi}

\end{document}